\RequirePackage{fix-cm}
\documentclass[11pt]{article}
\usepackage{graphicx}
\let\proof\relax
\let\endproof\relax

\usepackage{pbox}
\usepackage{setspace,graphicx,epstopdf,amsmath,amsfonts, amssymb, amsthm, xfrac, colortbl, tabu,bm}
\usepackage{marginnote,datetime,enumitem,rotating,fancyvrb}
\usepackage{float}
\usepackage[caption = false]{subfig}
\usepackage{tabularx}
\usepackage{changepage}
\usepackage{placeins}
\usepackage{booktabs}
\usdate
\usepackage{indentfirst} 
\usepackage{geometry}
\usepackage{chngcntr}
\usepackage{multirow}
\usepackage{xcolor}
\usepackage{hyperref}
\usepackage[justification=centering, singlelinecheck=false,labelfont=bf]{caption}

\usepackage[bottom]{footmisc}
\usepackage{pgfplots, pgfplotstable}
\pgfplotsset{compat=1.9}
\usetikzlibrary{spy}
\usetikzlibrary{backgrounds}
\usetikzlibrary{decorations}

\pgfdeclarelayer{background layer}
\pgfsetlayers{background layer,main}

\newcounter{definition}
\newcommand{\definition}[1]{
  \refstepcounter{definition}
  \noindent{\textbf{Definition \thedefinition.#1 }}}

\newcounter{proposition}
\newcommand{\proposition}[1]{
  \refstepcounter{proposition}
  \noindent{\textbf{Proposition \theproposition.#1 }}}

\newcounter{condition}

  \newcounter{corollary}
\newcommand{\corollary}[1]{
 \refstepcounter{corollary}
  \noindent{\textbf{Corollary \thecorollary.#1 }}}

  \newcounter{theorem}
\newcommand{\theorem}[1]{
  \refstepcounter{theorem}
  \noindent{\textbf{Theorem \thetheorem.#1 }}}

  \newcounter{lemma}
\newcommand{\lemma}[1]{
  \refstepcounter{lemma}
  \noindent{\textbf{Lemma \thelemma.#1 }}}

  \newcounter{remark}
\newcommand{\remark}[1]{
  \refstepcounter{remark}
  \noindent{\textbf{Remark \theremark#1: }}}

\newfloat{algorithm}{H}{}
\floatname{algorithm}{Algorithm}
\floatstyle{plaintop}
\restylefloat{algorithm}

\makeatletter
\long\def\ifnodedefined#1#2#3{%
    \@ifundefined{pgf@sh@ns@#1}{#3}{#2}%
}

\pgfplotsset{
    discontinuous/.style={
    scatter,
    scatter/@pre marker code/.code={
        \ifnodedefined{marker}{
            \pgfpointdiff{\pgfpointanchor{marker}{center}}%
             {\pgfpoint{0}{0}}%
             \ifdim\pgf@y>0pt
                \tikzset{options/.style={mark=*, fill=white}}
                \draw [densely dashed, line width=0.5pt]
                (marker-|0,0) -- (0,0);
                \draw[options] plot coordinates {(marker-|0,0)};
             \else
                \tikzset{options/.style={mark=none}}
             \fi
        }{
            \tikzset{options/.style={mark=none}}
        }
        \coordinate (marker) at (0,0);
        \begin{scope}[options, fill=black]
    },
    scatter/@post marker code/.code={\end{scope}}
    }
}

\makeatother

\title{New bounds for truthful scheduling on two unrelated selfish machines
}
\author{  Olga Kuryatnikova  \thanks{ \href{mailto:kuryatnikova@gmail.com}{kuryatnikova@gmail.com} (corresponding author),
ORCID:  0000-0001-8460-7296} \hspace{0.1cm} and Juan C. Vera \thanks{\href{mailto:j.c.veralizcano@uvt.nl}{j.c.veralizcano@uvt.nl}} 
 \and 
\small{ Department of Econometrics and Operations Research, Tilburg University, 5037 AB Tilburg, Netherlands}
}

\date{}

\begin{document}
\maketitle


%

\begin{abstract}
We consider the minimum makespan problem for $n$ tasks and two unrelated parallel selfish machines. Let $R_n$ be the best approximation ratio of randomized monotone scale-free algorithms. This class contains the most efficient algorithms known  for truthful scheduling on two machines. We propose a new $Min-Max$ formulation for $R_n$, as well as upper and lower bounds on $R_n$ based on this formulation. For the lower bound, we exploit pointwise approximations of cumulative distribution functions (CDFs). For the upper bound, we construct randomized algorithms using distributions with piecewise rational CDFs. Our method improves upon the existing bounds on $R_n$ for small $n$. In particular, we obtain almost tight bounds for $n=2$ showing that $|R_2-1.505996|<10^{-6}$.
\end{abstract}
\noindent Keywords: minimax optimization,  truthful scheduling, approximation, piecewise functions, algorithmic mechanism design

\section{Introduction and main results} \label{intro}
Scheduling on unrelated parallel machines is a classical problem in discrete optimization. In this problem one has to allocate $n$ independent, indivisible tasks to $m$ simultaneously working unrelated machines. The goal is to minimize the time to complete all the tasks. This time is called the makespan, and the scheduling problem is called the minimum makespan problem. Lenstra et al. \cite{Lenstra} proved that the problem is NP-complete and that a polynomial-time algorithm cannot achieve an approximation ratio less than $\frac{3}{2}$ unless $P=NP$. \par
We restrict ourselves to the case of $m=2$ machines. For this case there is a linear-time algorithm by Potts \cite{Potts} and a polynomial-time algorithm by Shchepin and Vakhania \cite{Shchepin} which provide $\frac{3}{2}$-approximations. Both algorithms use linear programming (LP) relaxations of integer programs and rounding techniques. 

We are interested in  the minimum makespan problem in the setting of algorithmic mechanism design. In this setting, every machine belongs to a rational agent who requires payments for performing tasks and aims to maximize his or her utility. Nisan and Ronen~\cite{Nisan} introduced this approach to model interactions on the Internet, such as routing and information load balancing. The minimum makespan problem is one of many optimization problems considered in algorithmic mechanism design. These include, among others, combinatorial auctions (see, e.g., ~\cite{CorrRound},~\cite{Dobzinski} and references therein) and graph theoretic problems, such as the shortest paths tree~\cite{ShortPath} and the maximum matching problem~\cite{Match}.

To solve the minimum makespan problem in algorithmic mechanism design, one can use an \emph{allocation mechanism}. An allocation mechanism consists of two algorithms: one allocates tasks to machines, and the other allocates payments to agents (the machines' owners). The goal of the mechanism is to choose a  task allocation algorithm that minimizes the makespan and a payment allocation algorithm that motivates the agents to act according to the wishes of the algorithm designer. For instance, it is desirable that the agents reveal their correct task processing times to the designer. We consider direct revelation mechanisms. These mechanisms collect the information about the processing times from each agent and allocate the tasks and payments based on this information according to a policy known to the agents in advance. To maximize their utilities, the agents can lie about processing times of tasks on their machines. As a result, direct revelation mechanisms may be hard to implement correctly.

To motivate the agents to tell the right processing times, one can use a \emph{truthful} direct revelation mechanism. With such mechanisms, telling the truth becomes a dominant strategy for each agent regardless of what the other agents do. This property guarantees that the processing times used to construct the mechanism are correct. There is a vast literature on truthful mechanisms~\cite{Sqrt2,Mualem,Nisan,Saks}.  Not all task allocation algorithms can be used in truthful mechanisms. For instance, there is no known truthful mechanism for the polynomial-time algorithms by Potts~\cite{Potts} and Shchepin and Vakhania \cite{Shchepin}. Finding the best approximation ratio for truthful scheduling on unrelated machines is one of the hardest fundamental questions in mechanism design.

Saks and Yu~\cite{Saks} showed that a task allocation algorithm can be used in a truthful mechanism if and only if the algorithm is \emph{monotone}. Intuitively, a task allocation algorithm is monotone if it assigns a higher load to a machine as long as the processing times on this machine decrease (see Section~\ref{sec:prelimSched} for the formal definition of monotonicity). In this paper we focus on monotone task allocation algorithms and do not consider the allocation of payments.

Nisan and Ronen~\cite{Nisan} show that no deterministic monotone algorithm can achieve an approximation ratio smaller than two, but randomized algorithms can do better in expectation. From here on we say that a randomized allocation algorithm has  a given property, e.g., monotonicity, if this property holds with probability one according to the distribution of the random bits of the algorithm. Randomized task allocation algorithms that are monotone in this sense give rise to universally truthful mechanisms considered in this paper.

Next, a deterministic algorithm is \emph{task-independent} when the allocation of any task does not change as long as the processing times of this task stay fixed.  Every deterministic monotone allocation  algorithm on two machines with a finite approximation ratio is task-independent (Dobzinski and Sundararajan~\cite{Dobzinski}). Therefore, if a given randomized algorithm has a finite expected approximation ratio, this algorithm is task-independent with probability one. Hence we can restrict ourselves without loss of generality to monotone and task-independent randomized algorithms to find the best truthful approximation ratio for two machines.

Finally, we restrict our attention to \emph{scale-free} algorithms. An algorithm is scale-free if scaling all processing times by some positive number does not influence the output. Following Lu~\cite{Lu}, we note that for $m = 2$, scale-freeness and allocation independence imply that the allocation of each task depends only on the ratio of this  task's processing times, which simplifies the analysis. Scale-free algorithms are widely used in the literature, and the latest most efficient algorithms for truthful scheduling on two machines by Chen et al.~\cite{Chen}, Lu~\cite{Lu} or Lu and Yu~\cite{LuYu1} are scale-free.  In the sequel we work with monotone, task-independent, scale-free (denoted by MIS) task allocation algorithms. These algorithms provide good upper bounds on approximation ratios in scheduling \cite{Chen,Lu,LuYu1,LuYu2,Nisan}.  Lu and Yu~\cite{Lu,LuYu1,LuYu2} present a way to construct a payment allocation procedure for MIS algorithms which results in truthful allocation mechanisms.

Denote by $R_n$ the best worst-case expected approximation ratio of randomized MIS algorithms for the makespan minimization on two machines with $n$ tasks. For simplicity, in the rest of the paper we call $R_n$ the best approximation ratio. Our approach is to formulate a mathematical optimization problem for $R_n$. This approach was not common until recently when several successful truthful or truthful in expectation mechanisms have been constructed using linear or nonlinear programs \cite{CorrRound,Chen,TruthLP,Lavi}. This paper continues the trend to combine optimization with mechanism design and has the following contributions:
\begin{enumerate}
\itemsep0em
\item A $Min-Max$ formulation for $R_n$, see \eqref{main_result} in Corollary~\ref{cor0}. 
\item A unified approach to construct upper and lower bounds on $R_n$, see Section~\ref{bounds}.  \\
In formulation \eqref{main_result}, the outer minimization is over multivariate cumulative distribution functions (CDFs) and the inner maximization is over the positive orthant in two dimensions. This problem is in general not tractable, so we build bounds on the optimal value. The lower bounds are the result of restricting the inner maximization to a finite subset of the positive orthant. To obtain the upper bounds, we restrict the outer minimization to the set of piecewise constant CDFs.
This is a general approach which could work for any $Min-Max$ problem that requires optimizing over a set of functions, not necessarily CDFs.
\item New upper and lower bounds on $R_n$ for $n \in \{2,3,4\}$ and the task allocation algorithms corresponding to the given upper bounds (see Table~\ref{tab1}). The resulting upper bounds are currently the best for all monotone algorithms (not only MIS) on two machines.
\begin{table}[H]
\centering
\caption{Bounds on $R_n$}
\label{tab1}
\begin{tabular}{|c|c|c|c|c|}
\hline
\multirow{2}{*}{$n$} & \multicolumn{2}{c|}{Lower bound}                                                                                   & \multicolumn{2}{c|}{Upper bound}                                                                                         \\ \cline{2-5}
                     & existing                                                                                          & \textbf{new} & existing                                                                                                & \textbf{new} \\ \hline
2                    & \multirow{3}{*}{\begin{tabular}[c]{@{}c@{}}1.505949\\  \scriptsize \cite{Lu} \end{tabular}} & 1.5059953       & \begin{tabular}[c]{@{}c@{}}1.5068\\  \scriptsize \cite{Chen}  \end{tabular}                  & 1.5059964       \\ \cline{1-1} \cline{3-5}
3                    &                                                                                                   & 1.5076         & \multirow{2}{*}{\begin{tabular}[c]{@{}c@{}}1.5861\\ \scriptsize \cite{Chen}  \end{tabular}} & 1.5238         \\ \cline{1-1} \cline{3-3} \cline{5-5}
4                    &                                                                                                   & 1.5195         &                                                                                                         & 1.5628         \\ \hline
\end{tabular}
\end{table}
\item Almost tight bounds on $R_2$ (see Table \ref{tab1}). \\
For $n=2$ tasks the initial problem \eqref{main_result} simplifies to problem \eqref{second_main_result}, where the outer minimization is over univariate CDFs. We use piecewise rational CDFs to obtain the upper and lower bounds with a gap not larger than $10^{-6}$.
\end{enumerate}

Notice that another, less restrictive way to define randomized algorithms would be to say that the properties hold in expectation over the random bits of the algorithm. Monotone in expectation task allocation algorithms can be used in truthful in expectation mechanisms. For certain classes of problems (e.g., for combinatorial auctions), one can convert LP relaxations with rounding into truthful in expectation mechanisms, see Azar et al.~\cite{CorrRound}, Elbassioni et al.~\cite{TruthLP} or  Lavi and Swamy~\cite{Lavi}. Truthful in expectation mechanisms could perform better in expectation than the universally truthful ones. In this paper we do not analyze the former type of mechanisms. We refer the reader to Auletta et al.~\cite{SinBits} and Lu and Yu~\cite{LuYu2} for more information on truthfulness in expectation.

The outline of the paper is as follows. In Section 2 we provide more details about randomized MIS algorithms for two machines, describe results from earlier research, and formulate a mathematical optimization problem for $R_n$. In Section 3 we exploit the symmetry of this problem to analyze the performance of MIS algorithms and to obtain our $Min-Max$ formulation \eqref{main_result} for $R_n$. In Section 4 we construct and compute bounds on the optimal value of the $Min-Max$ problem for several small~$n$. In Section 5 we analyze the case with two tasks in more detail to improve the bounds for this case. Section 6 concludes the paper. In Section 7 we provide the omitted proofs. All computations are done in MATLAB R2017a on a computer with the processor  Intel\textsuperscript{\textregistered} Core\textsuperscript{\tiny{TM}} i5-3210M CPU @ 2.5 GHz and 7.7 GiB of RAM. To solve linear programs, we use IBM ILOG CPLEX 12.6.0 solver. 

\section{Preliminaries} \label{sec:prelimSched}

Unless otherwise specified, lower-case letters denote numbers, bold lower-case letters denote vectors, and capital letters denote matrices. For a given positive number $n$,  let $\mathbb{R}^n$ be the set of real vectors with $n$ entries. The notations $\mathbb{R}^n_+$ and $\mathbb{R}^n_{++}$ refer to nonnegative and strictly positive real-valued vectors, respectively. For a given positive number $m$, we use $[m]$ to denote the set $\{1,\dots,m\}$. 
Unless otherwise stated, we use parentheses to denote vectors and brackets to denote intervals, e.g., $(x_1,x_2)$ is a vector while $[x_1,x_2]$ and $[x_1,x_2)$ are intervals.

The input into the minimum makespan problem with $m$ machines and $n$ tasks is a matrix of processing times $T=(T_{ij}), \ i\in [m], \ j\in [n]$. We describe the solution to the problem by a task allocation matrix  $X \in \{0,1\}^{m{\times}n}$, such that $X_{ij}=1$ if task $j$ is processed on machine $i$ and $X_{ij}=0$ otherwise. Now, for given $X$ and $T$, we define the makespan  $M_i$ of machine $i$ and the overall makespan $M$. 
\begin{align}
M_i(X,T):=\sum_{j=1}^n X_{ij}T_{ij}, \ \ M(X,T):=\max_{i \in [m]} M_i(X,T), \label{def:Mx} 
\end{align}
The optimal makespan for $T$ is
\begin{align}
M^*(T):= \min_{X \in \{0,1\}^{m{\times}n}} M(X,T). \label{opt_m}
\end{align}
For an allocation algorithm  $\mathcal{A}$ and an input matrix $T$, $X^{\mathcal{A},T}\in \{0,1\}^{m{\times}n}$ denotes the output of $\mathcal{A}$ on $T$. If $\mathcal{A}$ is randomized, $X^{\mathcal{A},T}$ is a random matrix,  and we use  $M(X^{\mathcal{A},T},T)$  to denote the expected makespan.
The worst-case (expected) approximation ratio of $\mathcal{A}$ equals the supremum of the ratio $ \frac{M(X^{\mathcal{A},T},T)}{M^*(T)}$ over all time matrices $T$. We refer the reader to Motwani and Raghavan~\cite{Motwani} for a comprehensive discussion on randomized algorithms.


\subsection{The best approximation ratio of randomized MIS algorithms} \label{mechanism}

According to the definition in Section~\ref{intro}, randomized MIS algorithms are monotone, task-independent and scale-free (MIS)  with probability one. Therefore, by fixing the random bits of a randomized MIS algorithm, we obtain a deterministic MIS algorithm with probability one. We provide next a formal  description of deterministic MIS algorithms.

A task allocation algorithm is monotone if for every two processing time matrices $T$ and $T'$ which differ only on machine $i$, $\ \sum_{j=1}^n (X^{\mathcal{A},T}_{ij}-X^{\mathcal{A},T'}_{ij})(T_{ij}-T'_{ij}) \le 0$ (see \cite{Sqrt2}). That is, the load of a machine increases as long as the processing times on this machine decrease.
An algorithm is task-independent if the allocation of a task depends only on its processing times. To be precise, for any two time matrices $T$ and $T'$  such that $T_{ij} =T'_{ij}$ for task $j$ and all $ i \in [m]$, the allocation of task  $j$ is identical,  i.e., $X^{\mathcal{A},T}_{ij}  = X^{\mathcal{A},T'}_{ij},  \text{ for all } i \in  [m]$.
An algorithm is scale-free if the multiplication of all processing times by the same positive number does not change the allocation. That is, for any $T \in \mathbb{R}^{m\times n}_{++}$ and $\lambda>0$, the output of the algorithm on the inputs $T$ and $\lambda T$ is identical. \\

\noindent Deterministic MIS algorithms for $m=2$ have been characterized by Lu~\cite{Lu}:

{\theorem[Lu~\cite{Lu}]\label{th-1} All deterministic MIS algorithms for scheduling on two unrelated machines are of the following form. For every task $j \in [n]$, assign a threshold
$z_j \in \mathbb{R}_{++}$ and one of the following two conditions: $T_{1j} < z_jT_{2j}$ or $T_{1j} \le  z_jT_{2j}$. The task goes to the first machine if and only if the corresponding condition is satisfied.}

Let $\mathcal{C}$ be the class of (randomized) algorithms which randomly assign a threshold $z_j$ and a condition $T_{1j} < z_jT_{2j}$ or $T_{1j} \le  z_jT_{2j}$ to each task $j$ and then proceed as given in Theorem~\ref{th-1} for the deterministic case. With probability one a randomized MIS algorithm is a MIS algorithm, and therefore of the form given by Theorem~\ref{th-1}. Hence, to find the best approximation ratio, it is enough to consider only algorithms in $\mathcal{C}$.  Next, we show that to find the best approximation ratio, we can restrict ourselves to a subclass of $\mathcal{C}$.

Let $\mathcal{P}_n$ be the family of Borel probability measures supported on the positive orthant, i.e.
 $\text{supp}(\mathbb{P}) \subseteq \mathbb{R}_{++}^n$ for all $\mathbb{P} \in \mathcal{P}_n$, where $\text{supp}(\mathbb{P})$ is the support of $\mathbb{P}$. We use $\mathbf{E}_{\mathbb{P}}[\ ]$ and $\mathbf{P}_{\mathbb{P}}[ \ ]$ to denote the expectation and the probability over the measure $\mathbb{P}$, respectively.  In the sequel we use the notions of a probability measure and the corresponding probability distribution interchangeably. For a $ \mathbb{P} \in \mathcal{P}_n$ we define algorithm $\mathcal{A}^{\mathbb{P}}$ as follows: \vspace{0.1cm}

\begin{algorithm}  \label{alg1} 
\caption{A monotone, task-independent, scale-free task allocation algorithm $\mathcal{A}^{\mathbb{P}}$ for two  machines}
{\normalsize{
\begin{center}
\begin{minipage}{.6\textwidth}
\begin{enumerate}
\itemsep0em
\item[]\hspace{-0.6 cm}  \textbf{Input}: processing time matrix $T=(T_{ij}) \in \mathbb{R}^{2\times n}_{++}$
\item[]\hspace{-0.6 cm} \textbf{Output}: allocation $X \in \{0, 1\}^{2\times n}$  \medskip
\item  Draw a vector of thresholds $(z_1 , z_2,\dots, z_n)$ from $\mathbb{P}$
\item For each task $j = 1, 2,\dots, n$ do
\item \hspace{0.3 cm} If $\frac{T_{1j}}{T_{2j}} < z_j$: $ \ X_{1j} \leftarrow 1, \ X_{2j} \leftarrow 0$
\item \hspace{0.3 cm}  Else: $ \ X_{1j} \leftarrow 0, \ X_{2j} \leftarrow 1$
\item Output $X$
\end{enumerate} 
\end{minipage}
\end{center}
}}
\end{algorithm}
Denote the family of all algorithms of the form above by $\mathcal{A}^{\mathcal{P}_n}$:
\[\mathcal{A}^{\mathcal{P}_n}:=\{\mathcal{A}^{\mathbb{P}} \ : \ \mathbb{P} \in \mathcal{P}_n \}.\]
Consider a measure $\mathbb{P} \in \mathcal{P}_n$ and the corresponding algorithm  $\mathcal{A}^{\mathbb{P}} \in \mathcal{A}^{\mathcal{P}_n}$. Let $X^{\mathbb{P},T}\in \{0,1\}^{2 \times n}$ be the randomized allocation produced by $\mathcal{A}^{\mathbb{P}}$ on time matrix $T$.  The expected makespan of $\mathcal{A}^{\mathbb{P}}$ on $T$ is
\begin{align}
M(\mathbb{P},T)=\mathbf{E}_{\mathbb{P}} \max \bigg \{ \sum_{j \in [n]} T_{2j}X_{2j}^{\mathbb{P},T},  \ \sum_{j \in [n]} \ T_{1j}X_{1j}^{\mathbb{P},T}\ \bigg \}. \label{def:M}
\end{align}
Recall that $M^*(T)$, defined in \eqref{opt_m}, denotes the optimal makespan for $T$. Let $R_n(\mathbb{P}, T)$ be the expected approximation ratio of $\mathcal{A}^{\mathbb{P}}$ on $T$ and $R_n(\mathbb{P})$ be the worst-case approximation ratio:\\
\begin{align}
R_n(\mathbb{P}, T)&=\frac{M(\mathbb{P}, T)}{M^*(T)}, \label{def:RnPT} \\
R_n(\mathbb{P}) &=\sup_{T \in \mathbb{R}_{++}^{2\times n}} \  R_n(\mathbb{P}, T) \label{def:RnP} 
\end{align}
\noindent It could happen that  for some $\mathbb{P}$ the ratio $ R_n(\mathbb{P},T)$ is unbounded in $T$. We do not consider these cases as we know that $R_n \le 1.5861$ (see Section~\ref{literature}). For ease of presentation, we work on $\overline{\mathbb{R}}_{+}=\mathbb{R}_{+} \cup \{\infty\}$ so that the supremum $\sup_{T \in \mathbb{R}_{++}^{2\times n}} \  R_n(\mathbb{P}, T)$ is always defined.

When a tie $\frac{T_{1j}}{T_{2j}} = z_j$ occurs for some $j \in [n]$, algorithms from the family $\mathcal{A}^{\mathcal{P}_n}$ send task $j$ to the second machine. In general, an algorithm in $\mathcal{C}$ could send the task to the first or the second machine. Next, we show that this behavior at the ties does not affect the worst-case performance:

{\theorem{}\label{thm:limit}  For a given number of tasks $n$, let $\mathbb{P} \in \mathcal{P}_n$ and define
\begin{align}
\mathcal{T}:=\big \{T\in  \mathbb{R}_{++}^{2\times n}:\mathbf{P}_{z \sim \mathbb{P}}\big[ z_j=\sfrac{T_{1j}}{T_{2j}}\big]=0 \text{ for all } j\in [n]\big \}. \label{def:niceTabs}
\end{align}
Let $R_n(\mathbb{P})$ be defined as in~\eqref{def:RnP}, then
\[R_n(\mathbb{P})=\sup_{T\in \mathcal{T}}R_n(\mathbb{P},T).\]
}
\noindent Theorem~\ref{thm:limit} is proven in Section~\ref{proof:thm:limit} and has the following implication:

\begin{corollary}{}\label{cor:limit}The best  approximation ratio over all randomized MIS algorithms is the best approximation ratio over all algorithms in $\mathcal{A}^{\mathcal{P}_n}$.
\end{corollary}

\noindent By Corollary~\ref{cor:limit}, the best approximation ratio over all randomized MIS algorithms is
 \begin{align}
R_n=\inf_{\mathbb{P} \in \mathcal{P}_n} \  R_n(\mathbb{P}).\label{p0}
\end{align}

Later in the paper we show that $R_n(\mathbb{P})$ is invariant under permutations of the tasks for every $\mathbb{P}\in \mathcal{P}_n$ (see Theorem \ref{th0}). Therefore, to compute $R_n$  using \eqref{p0}, we can restrict the optimization to the distributions $\mathbb{P}\in \mathcal{P}_n$ invariant under permutations of the random variables corresponding to the thresholds $(z_1,\dots,z_n)$ (see Theorem \ref{th1}).
For such distributions, $R_n(\mathbb{P})$ is determined by the worst-case performance of $\mathcal{A}^{\mathbb{P}}$ on each pair of the tasks. See Section~\ref{min_max_gamma}, and Proposition~\ref{prop4} in particular, for more detail. This means that for every $T\in \mathbb{R}^{2\times n}_{++}$, one can find the expected approximation ratio of $\mathcal{A}^{\mathbb{P}}$ by applying to all pairs of tasks the algorithm $\mathcal{A}^{\mathbb{P}_2}$, where $\mathbb{P}_2$ is the bivariate marginal distribution of $\mathbb{P}$ (by invariance, this distribution is the same for all pairs of thresholds). We use this property of family $\mathcal{A}^{\mathcal{P}_n}$ to construct problem \eqref{main_result} for $R_n$, which is one of the main results of this paper.


\subsection{Connection to the current knowledge on monotone algorithms} \label{literature}

The best approximation ratio for all monotone task allocation algorithms is unknown.  For deterministic algorithms with $n$ tasks and $m$ machines, when $n$ and $m$ tend to infinity, the ratio lies in the interval  $[1+\phi, \ m]$, where $\phi$ is the golden ratio.  The upper bound is due to Nisan and Ronen~\cite{Nisan}, and the lower bound is due to Koutsoupias and Vidali~\cite{Koutsoupias}. To compute this lower bound, the authors use a matrix of processing times where the numbers of rows and columns tend to infinity. If $n$ or $m$ is finite, the lower bound may be different. Koutsoupias and Vidali~\cite{Koutsoupias} present bounds for several finite time matrices as well.
For randomized algorithms, the best approximation ratio lies in the interval  $[2- \frac{1}{m}, \ 0.83685m]$.  The lower bound is due to  Mu'alem and Schapira~\cite{Mualem} who use Yao's minimax principle (for some details on this principle see, e.g., Motwani and Raghavan~\cite{Motwani}). The upper bound is due to Lu and Yu ~\cite{LuYu1}. The gap between the bounds grows with $m$, and the case with the smallest number of machines, $m=2$, has gained much attention in the literature.

For $m=2$, Nisan and Ronen~\cite{Nisan} have shown that the best approximation ratio of deterministic monotone algorithms is equal to $2$, for any finite $n$. The ratio for randomized monotone algorithms lies in the interval $[1.5,1.5861]$. Chen et al.~\cite{Chen} compute the upper bound using an algorithm from the family $\mathcal{A}^{\mathcal{P}_n}$ with independently distributed thresholds. The lower bound is the earlier mentioned bound by  Mu'alem and Schapira~\cite{Mualem}. There exist tighter lower bounds for certain cases. Lu~\cite{Lu} shows that algorithms from the family $\mathcal{A}^{\mathcal{P}_n}$ (and thus, by Corollary \ref{cor:limit}, all randomized MIS algorithms) cannot achieve a ratio better than $\frac{25}{16}$ ($=1.5625$) for sufficiently large $n$. Chen et al.~\cite{Chen} prove that an algorithm $\mathcal{A}^{\mathbb{P}}$ cannot do better than $1.5852$ when $\mathbb{P}$ is a product measure, i.e., when the thresholds are independent random variables. \par
The cases with $m=2$ and small $n>2$ are not well studied. Chen et al.~\cite{Chen} theoretically establish the upper bound $1.5861$ for any finite $n$. They also present numerical computations of upper bounds smaller than 1.5861 for some $n$. Finding these smaller bounds requires solving  non-convex optimization problems. However, the numerical method used by Chen et al.~\cite{Chen} does not guarantee global optimality. \par
The case with $m=2$ and $n=2$ is the simplest one, but even for this case the best approximation ratio is unknown. The ratio for algorithms from the family $\mathcal{A}^{\mathcal{P}_2}$  lies in the interval $[1.505949,1.5068]$. The upper bound is due to Chen et al.~\cite{Chen}, the lower bound is computed by Lu~\cite{Lu} using Yao's minimax principle. Notice that Lu~\cite{Lu} states that the lower bound is 1.506, but we repeated the calculations from this paper and obtained the number $1.505949$. Thus, when reporting results, we use this number as the currently best lower bound. We improve this bound and show that  $|R_2-1.505996|<10^{-6}$, in particular,  $R_2<1.506$.

From the description above, one can see that the existing bounds for truthful scheduling on unrelated machines are obtained using ad hoc procedures or (for the lower bounds) Yao's minimax principle. This paper develops a unified approach to construct upper and lower bounds on $R_n$ for any fixed $n$ and provides an alternative to Yao's minimax principle for the construction of lower bounds. As a result, we improve the bounds for truthful scheduling for $m=2$ and $n \in \{2,3,4\}$.

Next, we compare our approach to the existing methods for upper bounds in \cite{CorrRound,Chen,TruthLP,Lavi} that use optimization. Our method for upper bounds generalizes the approach by Chen et al.~\cite{Chen}. The  generalization considers a broader class of algorithms and provides stronger upper bounds for $n \le 4$. Our approach is fundamentally different from the  methods in \cite{CorrRound,TruthLP,Lavi}. To begin with, our method is suitable for  the minimum makespan problem on unrelated machines, while the methods in  \cite{CorrRound,TruthLP,Lavi} are not guaranteed to work for this problem. Next,  \cite{CorrRound,TruthLP,Lavi} use LP relaxations of the corresponding integer programs while we use the tools from continuous optimization to obtain possibly non-linear, but tractable approximations. Moreover, \cite{CorrRound,TruthLP,Lavi} consider truthful in expectation mechanisms only while we work with universally truthful ones.   

We obtain new bounds only for small $n\le 4$ because of the growing size of the lower bound optimization problems.  One can try to solve these problems in a more computationally efficient way using, for instance, the column generation technique (see, e.g., Gilmore and Gomory \cite{Gilmore}). As a result, one could  expect to obtain new bounds for $n>4$ with our approach since the solution to the upper bound problem is a relatively simple construction based on the solution to the lower bound problem, as described in Section~\ref{bounds}. However, improving the efficiency of the lower bound computation is out of the scope of the current paper.


\section{Using the symmetry of the problem to obtain a new formulation for the best approximation ratio} \label{symmetry}

In this section we exploit the fact that problem \eqref{p0} is invariant under permuting the tasks and the machines to simplify formulation~\eqref{p0} and obtain formulation \eqref{main_result} in Section \ref{min_max_gamma}.

\subsection{Using the symmetry of the problem} \label{subsec:sym}

 For a given number of tasks $n$, let $S_n$ be the symmetric, or permutation, group on $n$ elements.  Given a vector $\mathbf{z} \in \mathbb{R}^n$ and $\pi \in S_n$, the corresponding permutation of the  elements of $\mathbf{z}$ by $\pi$ is denoted by $\mathbf{z}\pi$. The group $S_n$  corresponds to column permutations in a given time matrix $T$. We define another group, $S_{inv}$, which corresponds to row permutations in $T$. $S_{inv}$ consists of the identity action $^{id}.$ and the action $^{inv}.$ which takes element-wise reciprocals of any vector  $\mathbf{x} =(x_1,\dots,x_n) \in \mathbb{R}_{++}^n$:
 \[^{id}\mathbf{x}= \mathbf{x}, \ \ ^{inv}\mathbf{x}= \big (\tfrac{1}{x_{1}},\tfrac{1}{x_{2}},\dots,\tfrac{1}{x_{n}} \big ).\]

Now, we define the action of $S_{inv}{\times}S_n$ on $\mathcal{P}_n$.
Given $\mathbb{P}\in \mathcal{P}_n$, $\gamma\in S_{inv}$, $\pi \in S_n$ and a random variable $\mathbf{z} \sim \mathbb{P}$, we consider the transformation $\mathbf{z} \rightarrow {^\gamma\mathbf{z}\pi}$. We define ${^\gamma\mathbb{P}\pi} \in \mathcal{P}_n$ as the distribution of ${^\gamma\mathbf{z}\pi}$. Next, we prove that problem  \eqref{p0} is convex and invariant under the action of $S_{inv}{\times}S_n$ on $\mathbb{P}$. As a result, to find the infimum in \eqref{p0}, it is enough to optimize over the distributions $\mathbb{P}$ invariant under the action of $S_{inv}{\times}S_n$. This approach is regularly used in convex programming, see Dobre and Vera~\cite{symCoposMatr}, Gatermann and Parrilo~\cite{symSDP} or de Klerk et al.~\cite{symStarAlg}.

Given distributions $\mathbb{P}_1,\dots,\mathbb{P}_k \in \mathcal{P}_n$, and weights $ \alpha_i\ge 0$  for all $ i \in [k]$ such that  $\sum_{i=1}^k\alpha_i=1$, we define the
 convex combination $\sum_{i=1}^k\alpha_i \mathbb{P}_i \in \mathcal{P}_n$ as the distribution where we draw from $P_i$, $i\in [n]$ with probability $\alpha_i$. 
The construction of $\sum_{i=1}^k\alpha_i \mathbb{P}_i$ and definitions~\eqref{def:M},\eqref{def:RnPT} imply that
\[R_n \bigg (\sum_{i=1}^k\alpha_i \mathbb{P}_i,T\bigg )=\sum_{i=1}^k\alpha_iR_n(\mathbb{P}_i,T).\]
Therefore, using~\eqref{def:RnP}, we have
\begin{align}
R_n\bigg (\sum_{i=1}^k\alpha_i \mathbb{P}_i\bigg) \le \sum_{i=1}^k\alpha_iR_n(\mathbb{P}_i), \label{ineq:conv}
\end{align}
that is,  $R_n(\mathbb{P})$ is convex in $\mathbb{P}$. Now, we show the invariance of $R_n(\mathbb{P})$ under the action of $S_{inv}{\times}S_n$.


{\theorem{} \label{th0} For any given number of tasks $n$,  $\mathbb{P} \in \mathcal{P}_n$, $\gamma\in S_{inv}$ and $\pi \in S_n$,
\[
 R_n(\mathbb{P})= R_n({^\gamma\mathbb{P}\pi}).\]
 }
\noindent The proof of Theorem~\ref{th0} is presented in Section~\ref{pr_th0}. 

\theorem{} \label{th1} For any given number of tasks $n$,
\begin{align}
R_n=\inf_{\mathbb{P} \in \mathcal{P}_n}  R_n(\mathbb{P}) \text{ such that } \mathbb{P} \text{ is invariant under the action of } S_{inv}{\times}S_n. \label{pr:inv}
\end{align}
\proof{}
As problem~\eqref{pr:inv} has a smaller feasibility set than problem~\eqref{p0}, the optimal value of problem~\eqref{pr:inv} is not smaller than $R_n$.  To prove the opposite inequality, we show that for any distribution $\mathbb{P} \in \mathcal{P}_n$ there is a distribution $\mathbb{Q} \in \mathcal{P}_n$ invariant under the action of $S_{inv}{\times}S_n$ such that $ R_n(\mathbb{Q})  \le R_n(\mathbb{P})$. Given $\mathbb{P} \in \mathcal{P}_n$,
 take $\alpha_i=\frac{1}{2(n!)}$ for $i \in [2(n!)]$ and consider the convex combination
\[\mathbb{Q}:=\frac{1}{2(n!)} \sum_{(\gamma,\pi) \in  S_{inv}{\times}S_n} {^{\gamma}\mathbb{P}\pi} . \]
By construction, $\mathbb{Q}$ has the required invariance property and
\begin{align*}
\hspace{-0.7cm} R_n(\mathbb{Q}) \overset{ \eqref{ineq:conv}}{\le} \tfrac{1}{2(n!)} \sum_{(\gamma,\pi) \in  S_{inv}{\times}S_n}  R_n({^{\gamma}\mathbb{P}\pi}) \overset{\text{Theorem }\ref{th0}}{=} \tfrac{1}{2(n!)} \sum_{(\gamma,\pi) \in  S_{inv}{\times}S_n}  R_n(\mathbb{P} )=R_n(\mathbb{P} ).
\end{align*}
\endproof
\subsection{New formulation for the best approximation ratio} \label{min_max_gamma}
From Theorem~\ref{th1}, problem \eqref{p0} is invariant under permuting the tasks and the machines. In the sequel we exploit the invariance under permuting the tasks only. First, this simplifies the presentation. Second, in our numerical computations using the invariance under the two types of permutations produced the same bounds as using invariance under task permutations only. 

Let $\mathcal{C}_n \subset \mathcal{P}_n$ be the family of probability measures invariant under the actions of $S_n$:
\begin{align}
\mathcal{C}_n=\{\mathbb{P}\in \mathcal{P}_n \ | \ \mathbb{P}={\mathbb{P}\pi}, \ \text{ for all } \ \pi \in S_n \}. \label{inv_distr}
\end{align}
In the rest of the paper we restrict the optimization to the distributions from $\mathcal{C}_n$. 
\corollary{} \label{cor:th1} For any given number of tasks $n$,
\begin{equation}\label{pr:inv1}
R_n=\inf_{\mathbb{P} \in \mathcal{C}_n}  R_n(\mathbb{P}).
\end{equation}
\proof{}
The Corollary follows from Theorem \ref{th1} and
\[\{\mathbb{P}\in \mathcal{P}_n \ | \ \mathbb{P}={^\gamma\mathbb{P}\pi}, \ \text{ for all } \ (\gamma,\pi) \in S_{inv}\times S_n \} \subset \mathcal{C}_n \subset \mathbb{P}.\]
\endproof

\noindent Proposition \ref{prop1} next is straightforward but crucial for our analysis.

{\proposition{} \label{prop1} Let $\mathbb{P} \in \mathcal{C}_n$. Then $\mathbb{P}$ has a cumulative distribution function (CDF) invariant under permutations of the variables.  Moreover, for $0<k<n$, all $k$-variate marginal distributions are identical. In particular, $\mathbb{P}$ is a joint distribution of $n$ identically distributed random variables.} 

By Proposition \ref{prop1}, if $\mathbb{P}\in \mathcal{C}_n$, then all univariate marginal distributions of $\mathbb{P}$ are identical and all bivariate marginal distributions of $\mathbb{P}$ are identical. Denote the corresponding univariate and bivariate CDFs by $F_{\mathbb{P}}$ and $H_{\mathbb{P}}$, respectively, and define
\begin{align}
\hspace{-0.7cm} \phi^{\mathbb{P}}(x,y){=}\hspace{0.1 cm} 1{+}y{-}\min\big\{ 1, \ 1{-}\sfrac{1}{x}+y\big\}F_{\mathbb{P}}(x){-}yF_{\mathbb{P}}(y){+}\min\big\{1{+}\sfrac{1}{x}, \ 1{+}y\big\}H_{\mathbb{P}}(x,y) \label{phi}
\end{align}

\noindent First, we present a result by Chen et al.~\cite{Chen}, which follows from Lu and Yu~\cite{LuYu2}

{\proposition[Chen et al.~\cite{Chen}] \label{prop4}  For any given number of tasks $n$, $\mathbb{P} \in \mathcal{C}_n$,  and $T\in \mathbb{R}_{++}^{2{\times}n}$,
\begin{align*}
R_n(\mathbb{P},T)  \le \max_{j,k \in [n]} \phi^{\mathbb{P}}\bigg(\frac{T_{1j}}{T_{2j}},\frac{T_{1k}}{T_{2k}} \bigg).
\end{align*}
}
Notice that this upper bound is defined by only two tasks out of $n$. Using Proposition~\ref{prop4}, we obtain the following formulation for $R_n(\mathbb{P})$:

{\theorem{} \label{th3} For any given number of tasks $n$, and $\mathbb{P} \in \mathcal{C}_n$, 
\begin{align*}
R_n(\mathbb{P})  = &
  \sup_{x,y \in \mathbb{R}_{++}}\phi^{\mathbb{P}}(x,y) .
\end{align*}
}
The proof of Theorem \ref{th3} is provided in Section \ref{pr_th3}.

{\remark{}\label{r1} Theorem \ref{th3} implies that the worst-case approximation ratio for $n$ tasks and $\mathbb{P} \in \mathcal{P}_n$ is the worst-case approximation ratio for two tasks and the bivariate marginal distribution of $\mathbb{P}$.} 

\noindent The next corollary is the main result of this section, and we use it throughout the rest of the paper. 

\corollary{} \label{cor0} For any given number of tasks $n$, 
\begin{align}
R_n= \inf_{\mathbb{P} \in \mathcal{C}_n} \ \sup_{x, y \in \mathbb{R}_{++}}  \phi^{\mathbb{P}}(x,y) . \label{main_result}
\end{align} 
\proof{}
The result follows from Corollary~\ref{cor:th1} and Theorem~\ref{th3}.
\endproof

\corollary{} \label{cor-1} $R_{n+1} \ge R_n$ for all $n\ge 2$.
\proof{}
The result follows from Corollary~\ref{cor0}.
\endproof


\section{Upper and lower bounds on the best approximation ratio} \label{bounds}
To find $R_n$ using \eqref{main_result}, one needs to optimize over a family of distributions, which is computationally intractable. Therefore we  construct upper and lower bounds on the optimal value of the problem. The idea is to restrict the attention to some subset of feasible distributions or some subset of $\mathbb{R}^{2}_{++}$, over which it is easier to solve problem \eqref{main_result}.
\begin{enumerate}
\item For the lower bound, we take a finite set  $\mathcal{S} \subset \mathbb{R}_{++}$ and find the supremum in \eqref{main_result} for $x,y \in \mathcal{S}$ only. A conventional approach to lower bounds is to propose several good-guess time matrices $T$, use these matrices to build a randomized instance of the minimum makespan problem,  and apply Yao's minimax principle. Our approach is different as we evaluate randomized algorithms on deterministic instances. 
\item For the upper bound, we find a good-guess distribution $\mathbb{P}$ and solve the inner maximization problem for this distribution. The distribution is built using the solution to the lower bound problem for $n \in \{2,3,4\}$. For $n=2$ we propose a more efficient approach in Section \ref{2t2m}.
\end{enumerate}

\subsection{Characterizing CDFs} 

To implement the ideas above, we have to optimize over distributions. For this purpose we represent a distribution via its CDF. To characterize CDFs, we follow Nelsen~\cite{Nelsen}. For $\mathbf{x}, \mathbf{y}\in \mathbb{R}^n$ such that $x_{i}\le y_i$ for all $i \in [n]$, we define the $n$-box 
$\mathcal{B}_{xy}: = \big [x_1,  y_1 \big ] \times \big [ x_2, y_2\big ] \times \dots \times \big [ x_n,  y_n \big ].
$
The set of vertices of $\mathcal{B}_{xy}$ is
$V_{xy}= \big \{x_1,  y_1 \big \} \times \big \{ x_2,  y_2\big \} \times \dots \times \big \{ x_n,  y_n \big \}.
$
The sign of vertex $\mathbf{b} \in V_{xy}$ is defined by
\[\text{sgn}(\mathbf{b}):=
 \begin{cases}
 \ \ 1, \ \text{if} \ b_i=x_{i} \ \text{for an even number of entries} \ i \\
-1, \ \text{if} \ b_i=x_{i} \ \text{for an odd number of entries} \ i.
  \end{cases}\]
Given a set $D \subseteq \mathbb{R}$, define $\overline{D}:=\big(D\cup \{0\}  \cup \{\infty\}\big)$. A function $G: \overline{\mathbb{R}}^n \rightarrow \mathbb{R}$ is called $n$-increasing on $\overline{D}^n$ when
\begin{align}
\sum_{\mathbf{b} \in V_{xy}} \text{sgn}(\mathbf{b}) G(\mathbf{b}) \ge 0,  \text{ for all } \ \mathbf{x}\le \mathbf{y} , \ \mathbf{x},\mathbf{y} \in \overline{D}^n \label{n_incr}
\end{align}
{\remark [Chapter 2.1 in Nelsen~\cite{Nelsen}]
For $n>1$, the fact that $G$ is $n$-increasing does not necessarily imply that $G$ is non-decreasing in each argument, and the other way round.
} 

\noindent The following family of functions captures the concept of  CDF.
\definition{\label{def:famG} Let $\mathcal{S} \subseteq \mathbb{R}_{++}$. $\mathcal{G}_n(\mathcal{S})$ is the family of functions $G: \overline{\mathcal{S}}^n \rightarrow [0,1]$ satisfying the conditions below.
\begin{enumerate}
\itemsep0em
\item $G$ is right continuous on $\overline{\mathcal{S}}^n$
\item $G$ is $n$-increasing on $\overline{\mathcal{S}}^n$
\item $G(\mathbf{z})=0$ for all $\mathbf{z}$ in  $\overline{\mathcal{S}}^n$ such that at least one of $z_i=0$
\item $G(\infty,\dots,\infty)=1$
\end{enumerate}}

{\proposition[Definition 2.10.8. in Nelsen~\cite{Nelsen}]\label{prop2}
A function $G: \overline{\mathbb{R}}_{++}^n \rightarrow [0,1]$ is a CDF of some $\mathbb{P}\in \mathcal{P}_n$ if and only if $G \in \mathcal{G}_n(\mathbb{R_{++}})$. } 

\subsection{Formulation of upper and lower bounds} \label{form_bounds}
 To construct the upper bound, we restrict the inner maximization in problem~\eqref{main_result} to a subset of $\mathbb{R}_{++}$. We do this using the next lemma.

\begin{figure}[H]
 \centering
\caption{$\mathcal{S}=\{\tfrac{1}{3},\tfrac{2}{3},1,\tfrac{3}{2},3\}$, $n=1$. The left plot: a function $g \in  \mathcal{G}_n(\mathcal{S})$. The right plot: a function $G \in \mathcal{G}_n(\mathbb{R}_{++})$. Notice that $G$ is a CDF and $g$ is the restriction of $G$ to $\overline{\mathcal{S}}^n$.}
\label{fig1}
\advance\leftskip-0.1cm
\begin{tikzpicture}[scale=0.48]
\begin{axis}[
    clip=false,
    jump mark left,
    ymin=0,ymax=1.2,
    xmin=0, xmax=4,
    x post scale=2,
    y post scale=1,
    xticklabel style={font=\LARGE},
    yticklabel style={font=\LARGE},
    xticklabels={$0$, \fontsize{18}{10}{$\sfrac{1}{3}$}, \fontsize{18}{10}{$\sfrac{2}{3}$}, $1$,  \fontsize{18}{10}{$\sfrac{3}{2}$},$3$},
    yticklabels={$0$,$0.25$,$0.45$,$0.55$,$0.75$,$0.9$},
    axis lines = middle,xlabel=$\bm{x}$,ylabel=$\bm{g}$,
    xlabel style={at={(axis description cs:1,-0.15)},font=\Large},
    ylabel style={at={(axis description cs:-0.06,1)},font=\Large},
    xtick={0,0.3333,0.6667,1,1.5,3},
    ytick={0,0.25,0.45,0.55,0.75,0.9},
    every axis plot/.style={very thick},
    discontinuous,
    table/create on use/cumulative distribution/.style={
        create col/expr={\pgfmathaccuma + \thisrow{f(x)}}
    },
    every tick/.style={
        black,
        semithick,
      }
]

\end{axis}
\draw[yscale=4.75, xscale=3.42]
node[fill=black, circle, inner sep=0pt,minimum size=3pt]{} (0.33333, 0.25)
node[fill=black, circle, inner sep=0pt,minimum size=3pt]{} (0.66667, 0.45)
node[fill=black, circle, inner sep=0pt,minimum size=3pt]{} (1, 0.55)
node[fill=black, circle, inner sep=0pt,minimum size=3pt]{} (1.5,0.75)
node[fill=black, circle, inner sep=0pt,minimum size=3pt]{} (3, 0.9)
node[fill=black, circle, inner sep=0pt,minimum size=3pt]{} (0, 0) ;
\end{tikzpicture}
\begin{tikzpicture}[scale=0.48]
\begin{axis}[
    clip=false,
    jump mark left,
    ymin=0,ymax=1.2,
    xmin=0, xmax=4,
    x post scale=2,
    y post scale=1,
    xticklabel style={font=\LARGE},
    yticklabel style={font=\LARGE},
    xticklabels={$0$,\fontsize{18}{10}{$\sfrac{1}{3}$}, \fontsize{18}{10}{$\sfrac{2}{3}$}, $1$,   \fontsize{18}{10}{$\sfrac{3}{2}$}, $3$, $a\hspace{0.05cm}{>}\hspace{0.05cm}3$},
    yticklabels={$0$,$0.25$,$0.45$,$0.55$,$0.75$,$0.9$,$1$},
    axis lines = middle,xlabel=$\bm{x}$,ylabel=$\bm{G}$,
    xlabel style={at={(axis description cs:1,-0.15)},font=\Large},
    ylabel style={at={(axis description cs:-0.07,1)},font=\Large},
    xtick={0,0.3333,0.6667,1,1.5,3,3.5},
    ytick={0,0.25,0.45,0.55,0.75,0.9,1},
    every axis plot/.style={very thick},
    discontinuous,
    table/create on use/cumulative distribution/.style={
        create col/expr={\pgfmathaccuma + \thisrow{f(x)}}
    },
    every tick/.style={
        black,
        semithick,
      }
]

\addplot [black] table [y=cumulative distribution, row sep=crcr]{%
 x f(x) \\
 0 0 \\
 0.33333 0.25 \\
 0.66667 0.2 \\
 1 0.1 \\
 1.5 0.2 \\
 3 0.15 \\
 3.5 0.10 \\
 4 0 \\
};
\end{axis}
\draw[yscale=4.75, xscale=3.42]
node[fill=black, circle, inner sep=0pt,minimum size=3pt]{} (0, 0);
\end{tikzpicture}
\end{figure}
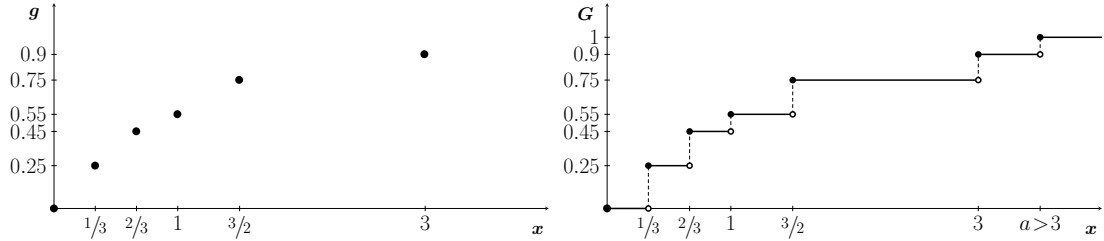

\lemma{}\label{l-1}  Let  $\mathcal{S} \subseteq \mathbb{R}_{++}$ be a \textbf{finite} set. Then $g \in  \mathcal{G}_n(\mathcal{S})$ if and only if there exists  $G \in \mathcal{G}_n(\mathbb{R}_{++})$ such that $g=G|_{\overline{\mathcal{S}}^n}$. That is, $g$ is a restriction of $G$ to $\overline{\mathcal{S}}^n$.
\proof{}
If there is $G \in \mathcal{G}_n(\mathbb{R}_{++})$ such that $g=G|_{\overline{\mathcal{S}}^n}$, then $g \in  \mathcal{G}_n(\mathcal{S})$ by definition of $\mathcal{G}_n(\mathbb{R_{++}})$. On the other hand, let $g \in  \mathcal{G}_n(\mathcal{S})$ and consider a number $a>\max \{s: s \in \mathcal{S}\}$. Let $\mathcal{S}_a=\mathcal{S} \cup \{a\}$  and define a new function $\hat{g}:\overline{\mathcal{S}_a}^n \rightarrow[0,1]$ such that $\hat{g}(\mathbf{z})=g(\mathbf{z})$ for $\mathbf{z} \in \overline{\mathcal{S}}^n$. For $\mathbf{z} \notin \overline{\mathcal{S}}^n$, construct a new vector $\mathbf{y}$ by replacing all occurrences of $a$ in $\mathbf{z}$ with $\infty$ and define $\hat{g}(\mathbf{z})=g(\mathbf{y})$. Consider the following piecewise constant function:
\begin{align}
G(z_1,\dots,z_n):= & \ \hat{g} \big (\max_{x \in \overline{\mathcal{S}_a}} \{x:  x\le z_1 \},\dots,\max_{x \in \overline{\mathcal{S}_a}} \{x:  x\le z_n \} \big ). \label{extention}
\end{align}
It is straightforward to show that $G \in \mathcal{G}_n(\mathbb{R}_{++})$ and $g=G|_{\overline{\mathcal{S}}^n}$. See Figure \ref{fig1} for an illustration of the case $n=1$.
\endproof

{\remark{}\label{r4}
The choice of the number $a$ is free and might influence our upper bound computations in Section~\ref{impl_bounds}.}\\

\noindent For a finite $\mathcal{S} \subset \mathbb{R}_{++}$ and  $g \in  \mathcal{G}_n(\mathcal{S})$, we define the restriction of the objective in problem~\eqref{main_result} to $\mathcal{S}$:
\begin{align}
\hspace{-0.7cm}\phi^{g}(x,y)=&\ 1+y-\min\big\{ 1, \ 1-\sfrac{1}{x}+y\big\}g(x,\infty,\dots,\infty)-yg(y,\infty,\dots,\infty) \label{phi_restr}\\[-4pt]
& \ +\min\big\{1+\sfrac{1}{x}, \ 1+y\big\}g(x,y,\infty,\dots,\infty)  \ \ \text{ for all } x,y \in \overline{\mathcal{S}}. \nonumber
\end{align}
By Lemma \ref{l-1},  $\phi^{g}=\phi^{\mathbb{P}}|_{\overline{\mathcal{S}}^2}$ for some $\mathbb{P} \in \mathcal{P}_n$.

\theorem{} \label{th4}
Given a number of tasks $n$, for any $\mathbb{P} \in \mathcal{C}_n$ and finite $\mathcal{S} \subset \mathbb{R}_{++}$, we have
\begin{equation}\label{lbub}
R_n(\mathbb{P} )\ge R_n \ge R_n(\mathcal{S}){:=}\inf_{g \in  \mathcal{G}_n(\mathcal{S})} \sup_{x,y\in \mathcal{S}} \hspace{-0.1cm}
\big  \{ \phi^{g}(x, y):  g(\mathbf{z})=g(\mathbf{z}\pi)  \text{ for all } \pi \in S_n,  \mathbf{z} \in \overline{\mathcal{S}}^n \big \}  
\end{equation}
\proof{}
The first inequality follows immediately from Corollary \ref{cor0}. Now, we prove the second inequality. Every $\mathbb{P} \in \mathcal{C}_n$ has a CDF $G_{\mathbb{P}} \in \mathcal{G}_n(\mathbb{R_{++}})$ invariant under permutations of the variables by Proposition~\ref{prop1}. At the same time, every such invariant $G \in \mathcal{G}_n(\mathbb{R_{++}})$  corresponds to some $\mathbb{P}_{G} \in \mathcal{C}_n$.  Combining this with Corollary~\ref{cor0}, we obtain
\begin{align*}
R_n= & \ \inf_{G \in \mathcal{G}_n(\mathbb{R_{++}})} \ \sup_{x, y \in \mathbb{R}_{++}} \{ \phi^{\mathbb{P}_G}(x,y)  :  G(\mathbf{z})=G(\mathbf{z}\pi)  \text{ for all } \pi \in S_n, \ \mathbf{z} \in \overline{\mathbb{R}}^n_{++} \} \\
\ge & \ \inf_{G \in \mathcal{G}_n(\mathbb{R_{++}})} \ \sup_{x, y \in \mathcal{S}} \{ \phi^{\mathbb{P}_G}(x,y)  :   G(\mathbf{z})=G(\mathbf{z}\pi)  \text{ for all } \pi \in S_n, \ \mathbf{z} \in \overline{\mathbb{R}}^n_{++} \}\\
=& \ \inf_{g \in \mathcal{G}_n(\mathcal{S})} \ \sup_{x, y \in \mathcal{S}} \{ \phi^{g}(x,y)  :  g(\mathbf{z})=g(\mathbf{z}\pi)  \text{ for all } \pi \in S_n, \ \mathbf{z} \in \overline{\mathcal{S}}^n \}.
\end{align*}
The last equality holds by Lemma \ref{l-1}. Notice that if $G \in \mathcal{G}_n(\mathbb{R_{++}})$ is invariant under permutations of the variables, then so is the $g:=G|_{\overline{\mathcal{S}}^n}$. On the other hand, if $g \in \mathcal{G}_n(\mathcal{S})$ is invariant under permutations of the variables, then so is the $G$ defined in \eqref{extention}.

\endproof

\subsection{Implementation and numerical results} \label{impl_bounds}

Let $S\subset \mathbb{R}_{++}$. To compute the lower bound $R_n(\mathcal{S})$ from \eqref{lbub}, we use the epigraph form of the optimization problem for $R_n(\mathcal{S})$:
\begin{align} \label{pr:lb_main}
R_n(\mathcal{S}) = &  \inf_{g \in  \mathcal{G}_n(\mathcal{S}), \ t \in \mathbb{R}} \ \ \ t &  \\
& \hspace{0.8cm} \text{s.t.}  \hspace{0.9cm} \phi^{g}(x, y) \le t & \text{for all } x,y\in \mathcal{S} \nonumber\\
& \hspace{2.2cm} g(\mathbf{z})=g(\mathbf{z}\pi) & \text{ for all } \pi \in S_n, \  \mathbf{z} \in \overline{\mathcal{S}}^n \nonumber\
\end{align}
The optimization variable in the problem above is $g$. This variable is a vector in $\mathbb{R}^{(|\mathcal{S}|+2)^n}$ which represents a function $g \in \mathcal{G}_n(\mathcal{S})$. We slightly abuse the notation and do not use a bold symbol for $g$ to underline that $g$ corresponds to a function with a finite support.  Family $\mathcal{G}_n(\mathcal{S})$ is an infinite family of functions $g$, and each of these functions is defined on a finite set $\overline{\mathcal{S}}^n$ with cardinality $(|\mathcal{S}|+2)^n$. For the purpose of optimization, this means that we consider all vectors $g\in \mathbb{R}^{(|\mathcal{S}|+2)^n}$ which satisfy the four conditions in Definition~\ref{def:famG} and the invariance property. All mentioned conditions are linear, and there are finitely many of them. Therefore the optimization over the infinite family of functions $\mathcal{G}_n(\mathcal{S})$ can be written as a finite LP. We use the invariance of $g$ (the second constraint) and Conditions 3-4 in Definition~\ref{def:famG} to reduce the number of variables problem~\eqref{pr:lb_main} (the size of $g$ as a vector). To ensure that $g$ corresponds to an $n$-increasing function as specified in \eqref{n_incr}, it is enough to consider $\mathbf{x}, \mathbf{y} \in \overline{\mathcal{S}}^n$ such that $x_i,y_i$ are consecutive points in $\mathcal{S}$ for all $i \in [n]$. This reduces the number of constraints in problem~\eqref{pr:lb_main}.

To compute the upper bound $R_n(\mathbb{P})$ using formulation \eqref{lbub}, we first construct a good-guess distribution $\mathbb{P}$. Given a set $\mathcal{S}$ and the solution $g$ to the lower bound problem \eqref{lbub} on $\mathcal{S}$, we use the distribution  $\mathbb{P}_{g}$ which corresponds to the CDF  \eqref{extention} based on $g$. To construct this CDF, we choose a number $a>\max \{s: s\in \mathcal{S}\}$, as explained in the proof of Lemma~\ref{l-1}. In the rest of this section we work with $\mathcal{S}_a=\mathcal{S} \cup \{a\}$. To solve \eqref{lbub} for $\mathbb{P}_{g}$, we define the following set of intervals:
\begin{align*}
\mathcal{I}_{\mathcal{S}}=& \ \big \{I_1,\dots,I_{|\mathcal{S}|+2}\big \} =\big \{[0,s_1),\dots,[s_{|\mathcal{S}|}, a), [a, \infty)\big\},
\end{align*}
\noindent This set of intervals covers $\mathbb{R}_{+}$, therefore by~\eqref{lbub}
\begin{align}
R_n(\mathbb{P}_{g}){=}&   \sup_{x,y\in \mathbb{R}_{++}}  \phi^{\mathbb{P}_{g}}(x, y) = \max_{I_i,I_j \in \mathcal{I}_{\mathcal{S}}} \bigg \{ \sup_{x \in I_i,y \in I_j} \phi^{\mathbb{P}_{g}}(x, y) \bigg\}. \label{ub_g}
\end{align}
We solve the inner maximization problem in \eqref{ub_g} for each pair $i,j \in [|\mathcal{S}|+2]$.
The expression for $\phi^{\mathbb{P}_{g}}$~\eqref{phi}  for the case $xy\ge 1$ is different from the case $xy<1$.   To simplify the computations when the line $xy=1$ crosses the rectangle $I_i{\times}I_j$, we restrict our attention to $\mathcal{S}$ of a particular type. Consider a collection of $k-1$ positive real numbers $r_1<r_2<\dots<r_{k-1}<1$, let
\begin{align}
\mathcal{S}_k=\{ s_1, \ s_2,\dots,s_{2k-1}  \}=\big\{r_1, \ r_2,\dots,r_{k-1}, \ 1, \ \tfrac{1}{r_{k-1}},\dots,\tfrac{1}{r_2},\ \tfrac{1}{r_1} \big \}. \label{points}
\end{align}
For any $a>\tfrac{1}{r_1}$, $\mathcal{S}_{k,a}=\mathcal{S}_k\cup \{a\}$ subdivides $\mathbb{R}_{+}$ into  $2k{+}1$ intervals $\mathcal{I}_{\mathcal{S}_k}$. First, consider a pair of intervals $I_i,I_j \in \mathcal{I}_{\mathcal{S}_k}$ such that $i \notin \{1,2k{+}1\}$ and $j
\notin \{1,2k{+}1\}$ (i.e., neither of them are the first or the last interval). Due to the choice of $\mathcal{S}_k$, the line $xy=1$ crosses the rectangle $I_i \times I_j$ if and only if $i+j=2k+1$. Denote the bivariate marginal CDF and the univariate marginal CDF of $\mathbb{P}_{g}$ by $H_{g}$ and $F_{g}$, respectively. Then
\begin{align}
\hspace{-0.9cm}\phi^{\mathbb{P}_g}(x,y){=}\begin{cases}\ 1{+}y{-}F_{g}(x){-}yF_{g}(y){+}\big(1{+}\sfrac{1}{x}\big ) H_g(x,y)
& i{+}j{\ge}2k{+}1, \ xy{\ge}1\\
\ 1{+}y\big(1{-}F_{g}(x){-}F_{g}(y){+}H_{g}(x,y)\big)\\
 \ \ {-}\big(1{-}\sfrac{1}{x}\big)F_{g}(x){+}H_{g}(x,y)& i{+}j {\le}2k+1, \ xy{<}1.
\end{cases} \label{phi_ub}
\end{align}
We construct the CDF of $\mathbb{P}_{g}$ using \eqref{extention}, therefore
\begin{align}
H_{g}(x,y)=& \ g(s_{i-1},s_{j-1},\infty,\dots,\infty) \ \text{ for all } (x,y) \in I_i\times I_j  \label{cdf_ub} \\
F_{g}(x)= & \ g(s_{i-1},\infty,\dots,\infty) \ \text{ for all } x \in I_i. \nonumber
\end{align}
That is, the marginal CDFs are constant on $I_i\times I_j$ and $I_i$, respectively.
As the range of a CDF is $[0,1]$, we conclude that for $x \in I_i,y \in I_j$, $\phi^{\mathbb{P}^{g}}(x, y)$ is non-increasing in $x$ and non-decreasing in $y$. The latter holds since for any $ \mathbb{P} \in \mathcal{P}_2$ invariant under $S_2$ and for any $x,y \in \mathbb{R}$,
\begin{align}
 & \hspace{-0.9cm}
 H_\mathbb{P}(x,y)=F_\mathbb{P}(x)+F_\mathbb{P}(y)-1+\mathbf{P}_{\mathbf{z} \sim \mathbb{P}}[z_1>x,z_2>y] \ge \max \big\{0,
 F_\mathbb{P}(x)+F_\mathbb{P}(y)-1\big\}. \label{h_ineq}
\end{align}
Hence the optimal value of the inner maximization problem in \eqref{ub_g} can be obtained by first substituting the CDFs \eqref{cdf_ub} into the function \eqref{phi_ub} and then substituting $x=s_{i-1}$, $y=s_{j}$.  Note that this optimum is not attained. For the case $i+j=2k+1$, i.e., when the line $xy=1$ crosses the rectangle $I_i \times I_j$, the result holds due to the choice of $\mathcal{S}_k$. \par
When $i \in \{1,2k{+}1\}$ or $j
\in \{1,2k{+}1\}$, the function $\phi^{\mathbb{P}_{g}}(x, y)$ simplifies, and we solve such cases separately. The solution approach resembles the one from the previous paragraph.

In numerical experiments we use uniform finite sets
\begin{align}
\mathcal{S}^u_k =\big\{\tfrac{1}{k}, \tfrac{2}{k}, \dots,\tfrac{k-1}{k}, \ 1, \
\tfrac{k}{k-1},\dots,k \big\}. \label{uniform_set}
\end{align}
Table \ref{tab2} shows the best obtained bounds and the $k$ we use to compute these bounds.
\begin{table}[H]
\centering
\caption{Bounds for the case of two machines and two, three, and four tasks based on Theorem \ref{th4}}
\label{tab2}
\begin{tabular}{|c|c|c|c|c|c|}
\hline
\multirow{2}{*}{$n$} & \multicolumn{2}{c|}{Lower bound}     & \multicolumn{2}{c|}{Upper bound}   & \multirow{2}{*}{$k$} \\ \cline{2-5}
                   & Current                   & New      & Current                 & New      &                    \\ \hline
2                  & \multirow{3}{*}{\begin{tabular}[c]{@{}c@{}}1.505949\\  \scriptsize~\cite{Lu} \end{tabular}} & 1.505980 & \begin{tabular}[c]{@{}c@{}}1.5068\\  \scriptsize~\cite{Chen} \end{tabular}                  & 1.5093 & 250                \\ \cline{1-1} \cline{3-6}
3                  &                           & 1.5076 & \multirow{2}{*}{\begin{tabular}[c]{@{}c@{}}1.5861\\  \scriptsize~\cite{Chen} \end{tabular}}                  & 1.5238  & 50                 \\ \cline{1-1} \cline{3-3} \cline{5-6}
4                  &                           & 1.5195   &                         & 1.5628   & 20                \\ \hline
\end{tabular}
\end{table}

\noindent We round the lower bounds $R_n(\mathcal{S}^u_k)$ down and the upper bounds $R_n(\mathbb{P}_{g})$ up. We verify all upper bounds with exact arithmetics using the MATLAB symbolic package and the following procedure. First, we obtain the optimal solution $g$ to problem~\eqref{pr:lb_main} and round the elements of the set $\mathcal{S}^u_k$ and the number $a$ to the $8^{th}$ digit. Next, we transform the rounded values into rational numbers and compute $R_n(\mathbb{P}_{g})$ as a rational number. By Lemma~\ref{l-1} and Theorem~\ref{th4}, the rounded $g$ provides the algorithm $\mathcal{A}^{\mathbb{P}_{g}}$ with the worst-case approximation ratio $R_n(\mathbb{P}_{g})$.

The upper bound for $n=2$ in Table~\ref{tab2} is worse than the best existing upper bound. We improve our result and obtain a new best exiting upper bound for $n=2$ in the next section.


\section{More precise bounds for two tasks} \label{2t2m}
In this section we analyze the case with $n=2$ tasks and $m=2$ machines in more detail. Now, to obtain an upper bound, we do not simply use some good-guess distribution as we did before, but we optimize over a subset of $\mathcal{C}_n$~\eqref{inv_distr}. Moreover, as a side result of this optimization, we obtain a non-uniform set $\mathcal{S}_k$ which produces a better lower bound than the one from  Table \ref{tab2} in the previous section. \par
Problem \eqref{main_result} simplifies for $n=2$. Given $F \in \mathcal{G}_1(\mathbb{R}_{++})$, define
\begin{align}
H(x,y):=\max \big\{0, \
F(x)+F(y)-1\big\}. \label{clayton}
\end{align}
$H$ is a copula, i.e., there is  $\mathbb{P}_{H,F} \in \mathcal{P}_2$ for which $H$ is its CDF and $F$ is its marginal CDF. See Nelsen~\cite{Nelsen} for the detailed description of copulas and their properties. Moreover, by construction $\mathbb{P}_{H,F} \in \mathcal{C}_2$. Therefore we can rewrite problem~\eqref{main_result} using univariate marginal CDF's.

\theorem{} \label{th5} Consider family $\mathcal{G}_1(\mathbb{R}_{++})$ from Definition~\ref{def:famG}.
\begin{align}
R_2  = & \  \inf_{F \in \mathcal{G}_1(\mathbb{R}_{++})} \ \sup_{x,y \in \mathbb{R}_{++}} \ 1+y-\min\big\{ 1, \ 1-\sfrac{1}{x}+y\big\}F(x)-yF(y) \label{best_r2} \\
& \hspace{0.7cm} +\min\big\{1+\sfrac{1}{x}, \ 1+y\big\}\max \big\{0, \
F(x)+F(y)-1\big\}. \nonumber
\end{align} 
\proof{}  Given $F\in \mathcal{G}_1(\mathbb{R}_{++})$, let
\begin{align*}
\phi^F(x,y) =& \ 1+y-\min\big\{ 1, \ 1-\sfrac{1}{x}+y\big\}F(x)-yF(y) \\
&  \hspace{0.2cm} +\min\big\{1+\sfrac{1}{x}, \ 1+y\big\}\max \big\{0, \
F(x)+F(y)-1\big\}.
\end{align*}
Consider any $\mathbb{P}\in \mathcal{C}_2$ with the univariate CDF $F_{\mathbb{P}}$. Define $\phi^{\mathbb{P}}(x,y)$ as in \eqref{phi}. From \eqref{h_ineq} for all $x,y \in \mathbb{R}_{++}$, 
\begin{align*}
\phi^{\mathbb{P}}(x,y)\ge \phi^{F_{\mathbb{P}}}(x,y).
\end{align*}
Thus
\begin{align*}
\hspace{-0.9cm} R_2\overset{\eqref{main_result}}{=} \inf_{\mathbb{P} \in \mathcal{C}_2} \ \sup_{x,y \in \mathbb{R}_{++}} \phi^{\mathbb{P}}(x,y) \ge \inf_{\mathbb{P} \in \mathcal{C}_2} \ \sup_{x,y \in \mathbb{R}_{++}}  \phi^{F_{\mathbb{P}}}(x,y) \ge \inf_{F\in \mathcal{G}_1(\mathbb{R}_{++})} \ \sup_{x,y \in \mathbb{R}_{++}}  \phi^F(x,y).
\end{align*}
On the other hand, for all $F \in \mathcal{G}_1(\mathbb{R}_{++})$ there is copula $H$ from \eqref{clayton} with the corresponding distribution $\mathbb{P}_{H,F} \in \mathcal{C}_2$. Hence
\begin{align*}
\hspace{-0.8cm} R_2 = \inf_{\mathbb{P} \in \mathcal{C}_2} \ \sup_{x,y \in \mathbb{R}_{++}} \phi^{\mathbb{P}}(x,y) & \le \inf_{F\in \mathcal{G}_1(\mathbb{R}_{++})} \ \sup_{x,y \in \mathbb{R}_{++}}  \phi^{\mathbb{P}_{H,F}}(x,y)\\
& =\inf_{F\in \mathcal{G}_1(\mathbb{R}_{++})} \ \sup_{x,y \in \mathbb{R}_{++}}  \phi^F(x,y).
\end{align*}
\endproof
{\remark{}\label{r3}  Nelsen~\cite{Nelsen} shows that for  $n>2$ the function
\[G(x_1,\dots,x_n)=\max \big\{0, \
\sum_{i=1}^n F(x_i)-n+1\big\}\]
 is not a CDF. We do not see other suitable $n$-variate CDFs which would have a bivariate margin $H$ from  \eqref{clayton}. As a result, the proof of Theorem~\ref{th5} fails for $n>2$.
}


\subsection{New upper bound for two tasks} \label{ub_2tasks}

To compute a new upper bound on $R_2$, we restrict the set of functions in problem \eqref{best_r2} to the family of piecewise rational univariate CDFs. We say that a function is piecewise rational if it can be written as a fraction where both the numerator and the denominator are polynomials. The domain of each CDF is subdivided into pieces by $\mathcal{S}_k$ defined in \eqref{points}. Given  $\mathcal{S}_k$, we introduce a collection of intervals:
\begin{align}
\hspace{-0.7cm}\mathcal{I}_{\mathcal{S}_k}=& \ \big \{I_1,I_2,\dots, I_k,I_{k+1},\dots,I_{2k-1},I_{2k}\big \} \nonumber \\
=& \ \big \{[0,r_1),[r_1,r_2),\dots,[r_{k-1},1),\big[1,\sfrac{1}{r_{k-1}}\big),\dots,\big[\sfrac{1}{r_{2}},\sfrac{1}{r_{1}}\big), \big[\sfrac{1}{r_1},+\infty\big) \big\}. \label{def:int2}
\end{align}
{\remark{ \label{r5}
We build the intervals using the points from $\mathcal{S}_k$ only. This is different from  Section~\ref{impl_bounds} where we use an additional number $a>\max{\mathcal{S}_k} $ to construct the intervals.} }

\noindent Given $\mathcal{I}_{\mathcal{S}_k}$ and a family of continuous functions $\mathcal{F}$, we consider CDFs which ``piecewisely" belong to $\mathcal{F}$.

{\definition{\label{def:famF}
$\mathcal{C_F}(\mathcal{S}_k)$ is a family of functions $F: \mathbb{R}_{++} \rightarrow [0,1]$ such that 
{\fontsize{12}{13}\selectfont{\begin{align}
F(x) =  \begin{cases}
f_1(\sfrac{1}{x}) & x \in I_1 \\
f_2(\sfrac{1}{x}) & x \in I_2  \\
\dots\\
f_k(\sfrac{1}{x})  & x \in I_k \\
1-f_k(x) & x \in I_{k+1} \\
\dots\\
1-f_{2}(x) & x \in I_{2k-1} \\
1-f_1(x) & x \in I_{2k},
\end{cases} \label{f0}
\end{align}}}
$f_1(x)=0$, $f_i(x)\in \mathcal{F}$, $f_i(x)$ is non-decreasing, $f_{k}(1)\le 0.5$, $0\le f_i(x_{i})\le f_{i+1}(x_{i})$  for all $i<k$.}
}
\noindent By construction, $F$ is a CDF and thus $\mathcal{C_F}(\mathcal{S}_k)\subset \mathcal{G}_1(\mathbb{R}_{++})$.
Restricting the minimization in problem~\eqref{best_r2} to $\mathcal{C}_{\mathcal{F}}(\mathcal{S}_k)$ provides an upper bound on $R_2$. We use the symmetry of $F$ to simplify the expression for this bound.

\proposition{} Define $\mathcal{S}_k$ as in~\eqref{points} and consider family $\mathcal{C}_{\mathcal{F}}(\mathcal{S}_k)$ from Definition~\ref{def:famF}. $R_2$ is not larger than \label{prop3}
\begin{align}
\hspace{-0.9cm}
R_2\big (\mathcal{C}_{\mathcal{F}}(\mathcal{S}_k)\big ) & =  \inf_{F \in \mathcal{C}_{\mathcal{F}}(\mathcal{S}_k)} \ \sup_{x,y \in \mathbb{R}_{++}} \ 1+y-\min\big\{ 1, \ 1-\sfrac{1}{x}+y\big\}F(x)-yF(y) \label{pr:ub2} \\
& \hspace{2cm} +\min\big\{1+\sfrac{1}{x}, \ 1+y\big\}\max \big\{0, \
F(x)+F(y)-1\big\}\nonumber \\
&=  \inf_{F \in \mathcal{C}_{\mathcal{F}}(\mathcal{S}_k)} \ \sup_{x,y \in \mathbb{R}_{++}, \ xy \ge 1} \ y-\sfrac{1}{x} +\big(1+\sfrac{1}{x}-y\big)F(y)+\sfrac{1}{x}F(x). \label{second_main_result}
\end{align}
\proof{}
Let $F \in \mathcal{C}_{\mathcal{F}}(\mathcal{S}_k)$. Problem~\eqref{pr:ub2} is a restriction of problem~\eqref{best_r2} to a smaller set of functions. Therefore the former problem defines an upper bound on $R_2$.

\noindent Next, we show the equality between \eqref{pr:ub2} and~\eqref{second_main_result}. Denote by $\phi^F$ the objective of problem~\eqref{pr:ub2}. We claim that for every $x,y\in \mathbb{R}_{++}$,
\[\phi^F(x,y) \le \sup_{x,y \in \mathbb{R}_{++}, \ xy\ge 1} \phi^F(x,y).\]
Let $\hat{x}>0, \hat{y}>0$. First, consider the case $\hat{x} \hat y\ge 1$. By construction of \eqref{f0}, $\hat x \hat y \ge 1$ implies
$ F(\hat x)+F(\hat y) \ge F(\hat x)+F\big(\sfrac{1}{\hat x}\big) \ge 1$. Then
\begin{align}
\phi^F(\hat{x},\hat y) = &\ 1+\hat y-F(\hat x)-\hat yF(\hat y)  +\big(1+\sfrac{1}{\hat x}\big) (F(\hat x)+F(\hat y)-1)& \nonumber \\
 =&\hspace{0.2 cm} \hat y-\sfrac{1}{\hat x}+\big(1+\sfrac{1}{\hat x}-\hat y\big)F(\hat y)+\sfrac{1}{\hat x} F(\hat x). \label{phi_init}
 \end{align}
Now, let $\hat{x}\hat{y} < 1$. By construction of $\mathcal{I}_{\mathcal{S}_k}$ in~\eqref{def:int2}, there are $I_i, I_j \in \mathcal{I}_{\mathcal{S}_k}$ such that $\hat{x} \in I_i, \ \hat{y}\in I_j$. The set $ \mathcal{S}_k$ is finite,  therefore there is a sequence $\{ (x_t,y_t)\}_{t=1}^\infty$ such that for all $t$ the following holds: $x_t \in I_i, \ y_t \in I_j$, $x_t \rightarrow \hat{x}^+, \ y_t \rightarrow \hat{y}^+$, $x_ty_t<1$ and $x_t,y_t \notin \mathcal{S}_k$. For all $x\in \mathbb{R}_{++}\setminus \mathcal{S}_k$ we have $F(x)+F\big(\sfrac{1}{x}\big)=1$. Hence for all $t$
 \begin{align}
\phi^F(x_t,y_t)= \ & \ 1+y_t-\big(1-\sfrac{1}{x_t}+y_t\big)F(x_t)-y_tF(y_t) \nonumber \\
 =\ & \ y_t\big(1-F(y_t)\big)+\big(1-\sfrac{1}{x_t}+y_t\big)\big(1-F(x_t)\big)+\sfrac{1}{x_t}-y_t& \nonumber \\
= \ &\ \sfrac{1}{x_t}-y_t+\big(1-\sfrac{1}{x_t}+y_t\big)F\big(\sfrac{1}{x_t}\big)+y_tF\big(\sfrac{1}{y_t}\big)\nonumber \\
 \overset{\eqref{phi_init}}{=}  & \ \phi^F\big(\sfrac{1}{y_t},\sfrac{1}{x_t}\big). \label{phi_temp}
\end{align}
Finally, $F$ is right continuous in $(\hat{x},\hat{y})$, and so is $\phi^F(\hat{x},\hat{y})$. Since $x_t\rightarrow \hat{x}^+,  y_t\rightarrow \hat{y}^+$,
\begin{align*}
\phi^F(\hat{x},\hat{y})=\lim_{t\rightarrow \infty} \phi^F(x_t,y_t)\overset{\eqref{phi_temp}}{=}\lim_{t\rightarrow \infty} \phi^F\big(\sfrac{1}{y_t},\sfrac{1}{x_t}\big) \le \sup_{x,y \in \mathbb{R}_{++}, \ xy > 1} \phi^F\big(x,y\big).
\end{align*}
The last inequality follows from $\big(\sfrac{1}{y_t}\big)\big(\sfrac{1}{x_t}\big) > 1$.
\endproof

\subsection{Implementing the new upper bound for two tasks} \label{impl_2}

In this subsection, for a given $\mathcal{S}_k$, we choose  $\mathcal{F}$ in Definition~\ref{def:famF} to be the family of linear univariate functions
\begin{align}
 \mathcal{F}:=\{ c^0+c^1x: \ c^0,c^1 \in \mathbb{R} \}. \label{f_def}
\end{align}
Then $\mathcal{C}_{\mathcal{F}}(\mathcal{S}_k)$ includes, in particular, the CDFs from~\cite{Chen,Lu,LuYu1,Nisan} where the authors use piecewise functions with domains subdivided into two, four, or six intervals. We observe that the upper bounds  are better when the domains are subdivided more times or when each piece has a more complex form than a constant function, i.e., when $c^1$ can be non-zero. We improve upon the existing upper bounds by using a larger number of pieces and  letting  $c^1$ be non-zero in each piece. Define
\begin{align}
\phi^F(x,y):= y-\sfrac{1}{x} +\big(1+\sfrac{1}{x}-y\big)F(y)+\sfrac{1}{x}F(x). \label{phi_f}
\end{align}
Let $\mathcal{X}:=\big \{ (x,y) \in \mathbb{R}^2_{++}{:} \ xy \ge 1 \big \}$. Consider two formulations for $ R_2\big (\mathcal{C}_{\mathcal{F}}(\mathcal{S}_k)\big )$ which are equivalent to problem~\eqref{second_main_result}
\begin{align}
\hspace{-1cm}
R_2\big (\mathcal{C}_{\mathcal{F}}(\mathcal{S}_k)\big )
=&  \  \inf_{F \in \mathcal{C}_{\mathcal{F}}(\mathcal{S}_k), \ t} \hspace{0.2cm} t \label{cut_plains_constr}   \\
&\hspace{0.8 cm}  \text{s.t.} \hspace{0.8 cm}\phi^F(x,y) \le t,  \  \text{for all } (x,y)\in \mathcal{X} \nonumber\\
=&  \  \inf_{F \in \mathcal{C}_{\mathcal{F}}(\mathcal{S}_k), \ t} \hspace{0.2cm}t \label{ub_2} \\
&\hspace{0.8 cm}  \text{s.t.} \hspace{0.5 cm} \sup_{x \in I_i,y\in I_j}\phi^F(x,y) \le t,  \ \text{for all } i+j \ge 2k+1.  \nonumber
\end{align}
The second equality follows from the equivalence of problems~\eqref{second_main_result} and~\eqref{pr:ub2}  since for $\mathcal{S}_k$  defined in \eqref{points} $xy \ge 1$ implies $x \in I_i, \ y \in I_j$ with $i+j \ge 2k+1$. We use problems \eqref{cut_plains_constr}  and \eqref{ub_2} to approximate $R_2\big (\mathcal{C}_{\mathcal{F}}(\mathcal{S}_k)\big )$ with high precision. Namely, we use relaxations of problem \eqref{cut_plains_constr} to compute lower bounds on $R_2\big (\mathcal{C}_{\mathcal{F}}(\mathcal{S}_k) \big )$, and we use feasible solutions to problem \eqref{ub_2} to compute upper bounds on $R_2\big (\mathcal{C}_{\mathcal{F}}(\mathcal{S}_k) \big )$. \par
For $F \in \mathcal{C}_{\mathcal{F}}(\mathcal{S}_k)$ satisfying~\eqref{f0} and~\eqref{f_def}, the variables in problem \eqref{cut_plains_constr} are $\big (t,\{c_i^0\}_{i=1}^k,\{c^1_i\}_{i=1}^k \big) $. This problem is LP with infinitely many constraints: each $(x,y) \in \mathcal{X}$ induces a linear constraint. Such problems can be well approximated using the cutting-plane approach introduced by Kelley~\cite{Kelley}. Namely, we start with a finite set $\mathcal{Y} \subset \mathcal{X}$ and restrict the set of constraints in \eqref{cut_plains_constr} to its finite subset generated by $(x,y) \in \mathcal{Y}$. As a result, we obtain a finite linear problem. Denote its optimal solution by  $(\underline{F}, \ \underline{t})$. Then $\underline{t}$ is a lower bound on $R_2\big (\mathcal{C}_{\mathcal{F}}(\mathcal{S}_k) \big )$.\par
Next, we substitute $\underline{F}$ in \eqref{ub_2} and find a feasible $t$.  We compute the supremum for each pair $i,j \in [2k]$ with $i+j\ge 2k+1$ using the Karush-Kuhn-Tucker (KKT) conditions. For each pair of intervals the inner maximization problem in~\eqref{ub_2} either is linearly constrained or satisfies the Mangasarian-Fromovitz constraint qualification. Therefore the optimum is among the KKT points, see, e.g., Section 3 in Eustaquio et al.~\cite{ConQual} for more details. All problems are simple and have similar structure. Therefore we do not write the KKT conditions explicitly, but consider all possible critical points from the first order conditions and from the boundary conditions. This set contains all the KKT points, and thus the optimal one. At the end we choose the point $(x,y)$ with the maximal value of $\phi^{\underline{F}}(x,y)$ among the critical points. See Section \ref{ap1} for more details.
Let $\overline{t}$ be the maximum of $\phi^{\underline{F}}(x,y)$ over all pairs of intervals. The solution $(\underline{F},\overline{t})$ is feasible for problem \eqref{ub_2}. Thus $\overline{t}$ is an upper bound on $R_2\big (\mathcal{C}_{\mathcal{F}}(\mathcal{S}_k) \big )$. Let $(x^*,y^*)$ be a point such that $\phi^{\underline{F}}(x^*,y^*)=\overline{t}$. If $|\overline{t}-\underline{t}|>10^{-8}$, we proceed from the beginning by restricting problem \eqref{cut_plains_constr} to the updated set $\mathcal{Y} \leftarrow \mathcal{Y}\cup \{(x^*,y^*)\}$. Otherwise we stop. \par
To obtain numerical results, we use uniform sets $\mathcal{S}^u_k$ of the form \eqref{uniform_set}. We work with family $\mathcal{C}_{\mathcal{F}}(\mathcal{S}^u_k)$ from Definition~\ref{def:famF}, where the underlying family of functions  $\mathcal{F}$ is defined in \eqref{f_def}. We initialize the cutting-plane procedure with $\mathcal{Y}=\big \{(x,y): x,y \in \mathcal{S}^u_k, \ xy\ge 1\big \}$. The best obtained upper bound is indicated in bold in Table~\ref{tab3}, it is stronger than the currently best upper bound $1.5068$.

\begin{table}[H]
\centering
\caption{ Upper bounds on the best approximation ratio for two tasks}
\label{tab3}
\begin{tabular}{|c|c|c|c|c|c|}
\hline
\multicolumn{1}{|l|}{$k$} & 5     & 10     & 16     & 50      & 100 \\ \hline
Upper bound on $R_2$                                & 1.5174 & 1.5096 & 1.5066 & 1.5060 &  \textbf{1.5059964}\\ \hline
\end{tabular}
\end{table}
\noindent We verify the upper bound $1.5059964$ using exact arithmetics in a similar way as we do it for the upper bounds in Table \ref{tab2} of Section \ref{impl_bounds}.


\subsection{New lower bound for two tasks}  \label{impr_lb}
The cutting-plane approach from Section~\ref{ub_2tasks} generates not only the upper bound with the corresponding CDF, but also the set of points $\mathcal{Y}$. Using $\mathcal{Y}$, we build a new set $\mathcal{S}^*_k$ of the form \eqref{points}, which is not uniform as in \eqref{uniform_set}. We consider all $(x,y) \in \mathcal{Y}$ involved in the binding constraints of problem \eqref{cut_plains_constr} at the last cutting-plane iteration. Next, we take the corresponding $x$, $y$ and their reciprocals, order ascending, round to the $8^{th}$ digit and obtain $\mathcal{S}^*_k$ with $k=82$.  For this set the lower bound $R_n(\mathcal{S}^*_k)$ from problem~\eqref{pr:lb_main} is $1.5059953$, which is stronger than our lower bound from Table~\ref{tab2}. As a result, the lower and upper bounds become very close to each other: $| R_2-1.505996| \le 10^{-6}$.


\section{Conclusion}\label{conclusion}

We consider randomized MIS algorithms to the minimum makespan problem on two unrelated parallel selfish machines. We propose a new $Min-Max$ formulation \eqref{main_result} to find $R_n$, the best approximation ratio  over randomized MIS algorithms. The minimization in  \eqref{main_result} goes over distributions and the maximization goes over $\mathbb{R}^2_{++}$. The problem is generally intractable, so we build upper and lower bounds on the optimal value. To obtain the lower bound, we solve the initial problem on a finite subset of $\mathbb{R}^2_{++}$. Using the resulting solution, we construct a piecewise constant cumulative distribution function (CDF) for which the worst-case performance is easy to estimate. In this way, we obtain the upper bound on $R_n$. We implement this approach and find new bounds for $n \in \{2,3,4\}$ tasks. \par
For $n=2$ the best CDF is a known function of univariate margins (copula). We parametrize these margins as piecewise rational functions of degree at most one. The resulting upper bound problem \eqref{second_main_result} is a linear semi-infinite problem. We solve it by the cutting-plane approach. This approach provides the upper bound $1.5059964$ and the CDF for which the algorithm achieves this bound. \par
As a side result of the cutting-plane approach, we obtain the lower bound $1.5059953$, so $|R_2-1.505996| \le 10^{-6}$. 
This work leaves open several questions for future research.  First, our approach could be made more numerically efficient to provide better bounds for $m=2$ machines. For example,  column generation could solve lower bound problem \eqref{lbub} on denser grids and for the larger number of tasks $n$. Parametrizing distributions of more than two variables could improve the results for the upper bound problem \eqref{lbub}. Second, we work with $m=2$, machines but there are algorithms for $m>2$ machines with similar properties, e.g., by Lu and Yu~\cite{LuYu2}. It would be interesting to see how our approach works in the case of more than two machines. Finally, the suggested piecewise and pointwise constructions could be suitable for other problems with optimization over low dimensional functions, including problems from algorithmic mechanism design. 


 \section{Proofs}\label{proofs}

\noindent In this section we use some additional notation. First, $\rho(T)$ denotes the vector of the processing time ratios
\begin{align}
\rho(T):=\big ( \tfrac{T_{11}}{T_{21}},\tfrac{T_{12}}{T_{22}},\dots,\tfrac{T_{1n}}{T_{2n}} \big ). \label{rho_t}
\end{align}
Second, for $\mathbf{z} \in \mathbb{R}^n_{++}$ we denote by $\mathcal{A}^{\mathbf{z}}$ the algorithm in $\mathcal{A}^{\mathcal{P}_n}$ with the thresholds fixed at $\mathbf{z}$. Finally, for $T\in \mathbb{R}_{++}^{2\times n}$  we let $X^{\mathbf{z},T}$ and $M(\mathbf{z},T)$ be the output and the makespan of $\mathcal{A}^{\mathbf{z}}$ on $T$, respectively.

\subsection{Proof of Theorem~\ref{thm:limit}} \label{proof:thm:limit}

\noindent We begin with a lemma.

\lemma{} \label{lem:limit} For a given number of tasks $n$, let $\mathbb{P} \in \mathcal{P}_n$. For every time matrix $\ T\in \mathbb{R}_{++}^{2\times n}$, there exists a sequence of time matrices $\{ T_k\}_{k>0} \subset \mathbb{R}_{++}^{2\times n}$ such that $\mathbf{P}_{z \sim \mathbb{P}}\big[ z_j=\sfrac{T_{k,1j}}{T_{k,2j}}\big]=0$ for all $j\in [n]$ and
 $k$, and
 \[R_n(\mathbb{P}, T)=\lim_{k \rightarrow \infty} R_n(\mathbb{P}, T_k).\]
\proof{}
Consider a sequence of nonnegative numbers  $\{\epsilon_k\}$ such that
\begingroup
 \setlength{\leftmargini}{18pt}
\begin{enumerate}
\item $\mathbf{P}_{z\sim \mathbb{P}}\big [z_j={\rho(T)_j}+\epsilon_k \big ]=0, \ \forall \ k \in \mathbb{N}, \ j \in [n]$
\item $\displaystyle \lim_{k \rightarrow \infty} \epsilon_k=0$
\end{enumerate}
\endgroup
\noindent A sequence with these properties exists since for every $j\in [n]$ the case $\mathbf{P}_{z \sim \mathbb{P}}\big[ z_j=a\big]>0$ is possible for countably many  $a\in \mathbb{R}_{++}$ only. Next, we build a sequence of time matrices $\{ T_k\}_{k>0}, \ T_k=(T_{k,{ij}})$:
\begin{align*}
T_{k,{1j}}=T_{1j}+\epsilon_k T_{2j} \ \text{ and } \ T_{k,{2j}}=T_{2j} \ \text{ for all } j\in [n].
\end{align*}
Notice that $T=\lim_{k \rightarrow \infty} T_k$. By adding $\epsilon_k T_{2j}$ to every task on the first machine, we ensure that $\mathbf{P}_{z \sim \mathbb{P}}\big[ z_j=\sfrac{T_{k,1j}}{T_{k,2j}}\big]=0$ for all $j\in [n]$ and $k$. For each $k$ and all $j \in [n], \ i \in \{1,2\}$, we have $T_{ij}\le T_{k,{ij}}$. So  $M^*(T) \le M^*(T_k) \le M(X^{*},T_k)$, where $X^{*}$ is the optimal allocation for $T$. $X^{*}$ is finite (binary, in particular) and $T=\lim_{k \rightarrow \infty} T_k$. Combining this with \eqref{def:Mx}, we see that $M(X^{*},T_k)$ tends to $ M^*(T)$ when $k$ tends to infinity. Therefore
\begin{align}
 \lim_{k \rightarrow \infty}M^*(T_k)=M^*(T). \label{opt_lim}
\end{align}
For every time matrix $T_k$ and task $j$, consider the event  $B_{k,j}:`` \rho(T)_j<z_j\le \rho(T_k)_j"$. Let $B_k=\bigcup_{j=1}^n B_{k,j} $ and let $B_k^c$ be the complement of $B_k$. When $B_k$ happens, outcomes of $\mathcal{A}^{\mathbf{z}}$ on $T$ and $T_k$ are different, otherwise they are the same. By construction of $\mathcal{A}^{\mathbf{z}}$, $M(\mathbf{z},T)$ is finite for any $T$. Hence
\begin{align}
\mathbf{E}_{\mathbf{z} \sim \mathbb{P}} \big[M(\mathbf{z},T_k) -M(\mathbf{z},T)  \ | \ B_k^c \big ]=0, \ \ \text{ for all } \
 k \in \mathbb{N}. \label{event}
\end{align}
For any $j \in [n]$ we have  $ \rho(T_k)_j \rightarrow \rho(T)_j^+$. Since the CDF of $\mathbb{P}$ is continuous from the right,
\begin{align*}
& \lim_{k \rightarrow \infty} \mathbf{P}_{\mathbf{z} \sim \mathbb{P}}\big [ z_j\le \rho(T_k)_j \big ]=\mathbf{P}_{z \sim \mathbb{P}}\big [z_j\le \rho(T)_j\big ],
\end{align*}
which implies
\begin{align}
 \lim_{k \rightarrow \infty} \mathbf{P}_{\mathbf{z} \sim \mathbb{P}} [B_k]= & \ \lim_{k \rightarrow \infty} \mathbf{P}_{\mathbf{z} \sim \mathbb{P}}\bigg [\bigcup_{j=1}^n B_{k,j}\bigg ] \le \sum_{j=1}^n \lim_{k \rightarrow \infty} \mathbf{P}_{\mathbf{z} \sim \mathbb{P}}\big [B_{k,j}\big ] \nonumber \\
= & \ \sum_{j=1}^n \lim_{k \rightarrow \infty} \bigg ( \mathbf{P}_{\mathbf{z} \sim \mathbb{P}}\big [ z_j \le \rho(T_k)_j\big ]-\mathbf{P}_{\mathbf{z} \sim \mathbb{P}}\big [ z_j \le \rho(T)_j\big ] \bigg) =0.\label{zero_prob}
\end{align}
Finally, for any $k \in \mathbb{N}$
\begin{align*}
M(\mathbb{P},T_k)  - M(\mathbb{P},T)= & \ \mathbf{E}_{\mathbf{z} \sim \mathbb{P}} \big[ M(\mathbf{z},T_k)-M(\mathbf{z},T) \ | \ B_k \big ]\mathbf{P}_{\mathbf{z} \sim \mathbb{P}} [B_k] \\
 & \ + \mathbf{E}_{\mathbf{z} \sim \mathbb{P}} \big[ M(\mathbf{z},T_k)-M(\mathbf{z},T) \ | \ B_k^c \big ]\big (1-\mathbf{P}_{\mathbf{z} \sim \mathbb{P}}[B_k] \big )\\
\overset{ \eqref{event}}{=} & \ \mathbf{E}_{\mathbf{z} \sim \mathbb{P}} \big[ M(\mathbf{z},T_k)-M(\mathbf{z},T) \ | \ B_k \big ]\mathbf{P}_{\mathbf{z} \sim \mathbb{P}} [B_k] \\
\le & \ |T|_1\mathbf{P}_{\mathbf{z} \sim \mathbb{P}} [B_k]
\end{align*}
Thus by \eqref{zero_prob},
\begin{align}
\lim_{k \rightarrow \infty}  M(\mathbb{P},T_k)= M(\mathbb{P},T), \label{lim_M}
\end{align}
and
\begin{align*}
R_n(\mathbb{P}, T)= \frac{M(\mathbb{P},T)}{M^*(T)}\overset{\eqref{opt_lim}, \eqref{lim_M}}{=}\lim_{k \rightarrow \infty} \frac{M(\mathbb{P},T_k)}{M^*(T_k)}=\lim_{k \rightarrow \infty} R_n(\mathbb{P},T_k).
\end{align*} 
\endproof

\proof[\textbf{Proof of Theorem~\ref{thm:limit}} ]
Recall that $\mathcal{T}$ is defined in~\eqref{def:niceTabs}. By Lemma~\ref{lem:limit}, for every $T\in  \mathbb{R}_{++}^{2\times n}$ there exists a sequence of time matrices $\{ T_k\}_{k>0} \subset \mathcal{T}$ such that
 \[R_n(\mathbb{P}, T)=\lim_{k \rightarrow \infty} R_n(\mathbb{P}, T_k)\le \sup_{T \in  \mathcal{T}}R_n(\mathbb{P}, T).\]
Hence \[R_n(\mathbb{P})=\sup_{T \in  \mathbb{R}_{++}^{2\times n}}R_n(\mathbb{P}, T) \le \sup_{T \in  \mathcal{T}}R_n(\mathbb{P}, T).\]
The opposite inequality holds since $\mathcal{T}\subset \mathbb{R}_{++}^{2\times n}$.

\endproof

\subsection{Proof of Theorem~\ref{th0}}\label{pr_th0}
\noindent Recall the definitions of $S_{inv}$ and $S_n$ from Section \ref{symmetry}. Let $\sigma_{id}$ be the identity action of $S_2$ and $\sigma_{swap}$ be the non-identity action of $S_2$. We say that that the group $S_2{\times}S_n$ acts on a matrix in $\mathbb{R}^{2\times n}$ by
 permuting the $2$ rows and the $n$ columns of this matrix. Namely, for $A\in \mathbb{R}^{2\times n}$ and $(\sigma,\pi) \in S_2{\times}S_n$, we define the action of $(\sigma,\pi) \in S_2{\times}S_n$ on $A$ by
  \[\sigma A \pi :=(A_{\sigma i,\pi j}).\]
We show that for any $T \in \mathbb{R}^{2\times n}_{++}$ the  optimal makespan  $M^*(T)$ is invariant under the actions of $S_2{\times}S_n$ on $T$. Moreover, the expected makespan $M(\mathbb{P},T)$ is invariant under the actions of $S_2{\times}S_n$ on $T$ and $S_{inv}{\times}S_n$ on $\mathbb{P}$. These two results together imply Theorem \ref{th0}.

{\lemma{} \label{l0} $M^*(T)=M^*(\sigma T \pi)$ for all $(\sigma,\pi)  \in S_2{\times}S_n, \ T \in \mathbb{R}_{++}^{2\times n}$}
 \begin{proof} Consider a time matrix $T$ and actions $\pi \in S_n, \ \sigma  \in S_2$.  Let $X^{*}=(X^{*}_{ij})$ be an optimal allocation matrix for $T$. Then
$\sigma T \pi=(T_{\sigma i, \pi j}), \ \ \sigma X^* \pi=(X^{*}_{\sigma i, \pi j})$.
This implies
 \begin{align*}
& \ M^*(\sigma T \pi ) \le  \ \max_{i} \big \{ \sum_{j \in [n]} T_{\sigma i, \pi j}X^{*}_{\sigma i, \pi j} \big \} = \max_{i} \big \{ \sum_{j \in [n]} T_{ij}X^{*}_{ij} \big \}=M^*(T).
\end{align*}
Analogously, for the time matrix $\sigma T \pi$ and actions $\pi^{-1} \in S_n, \ \sigma^{-1} \in S_2$, we obtain
$$M^*( T ) \le M^*(\sigma T \pi).$$
 \end{proof}

\proof[\textbf{Proof of Theorem~\ref{th0}}] Let $\mathbb{P}\in \mathcal{P}_n$, $\mathbf{z} \in \mathbb{R}^n_{++}$ and $T\in \mathbb{R}^{2\times n}_{++}$. Since we are interested in $R_n(\mathbb{P})$, by Lemma~\ref{lem:limit} we can assume without loss of generality that $T$ is such that $\mathbf{P}_{z \sim \mathbb{P}}\big[ z_j=\sfrac{T_{1j}}{T_{2j}}\big]=0$ for all $j\in [n]$.  Consider $(\gamma,\pi) \in S_{inv}{\times}S_n$ and $\mathbf{y}={^\gamma\mathbf{z}\pi}$.
Let $\sigma=\sigma_{id}$ if $\gamma=^{id}.$ and $\sigma=\sigma_{swap}$ if $\gamma=^{inv}.$  Then $\mathcal{A}^{\mathbf{z}}$ sends task $j$ to machine $i$ on  $T$ if and only if $\mathcal{A}^{\mathbf{y}}$ sends task $\pi j$ to machine $\sigma i$ on $\sigma T \pi$.  As a result, $T_{ij}X^{\mathbf{z},T}_{ij} =T_{\sigma i \pi  j}X^{\mathbf{y},\sigma T \pi}_{\sigma i \pi  j} $ for all $i,j$ and
\begin{align*}
M\big (\mathbf{z},T \big  )= & \ \max_{i} \bigg \{ \sum_{j \in [n]} T_{ij}X^{\mathbf{z},T}_{ij} \bigg \}=\max_{i} \bigg \{ \sum_{j \in [n]} T_{\sigma i \pi  j}X^{\mathbf{y},\sigma T \pi}_{\sigma i \pi  j}  \bigg \} \nonumber \\
= & \ \max_{i} \bigg \{ \sum_{j \in [n]} T_{ i j}X^{\mathbf{y},\sigma T \pi}_{i j} \bigg \}= M\big (\mathbf{y},\sigma T \pi \big  ).
\end{align*}
Therefore
\begin{align}
M(\mathbb{P},T)= \mathbf{E}_{\mathbf{z} \sim \mathbb{P}} M\big (\mathbf{z},T \big  )= & \ \mathbf{E}_{\mathbf{y} \sim {^{\gamma}\mathbb{P}\pi}}  M\big (\mathbf{y},\sigma T \pi \big  )=M({^{\gamma}\mathbb{P}\pi},\sigma T \pi), \label{problem}
 \end{align}
 Combining this with  Lemma \ref{l0}, obtain
\begin{align*}
 R_n(\mathbb{P}, T) \ = R_n({^{\gamma}\mathbb{P}\pi},\sigma T \pi).
\end{align*}

\noindent By Theorem~\ref{thm:limit},
\begin{align*}
R_n(\mathbb{P}) =\sup_{T \in  \mathcal{T}}R_n(\mathbb{P},T)=\sup_{T \in  \mathcal{T}}R_n({^{\gamma}\mathbb{P}\pi},\sigma T \pi)=R_n({^{\gamma}\mathbb{P}\pi}),
\end{align*}
where $\mathcal{T}$ is defined in~\eqref{def:niceTabs}.
\endproof


\subsection{Proof of Theorem \ref{th3}} \label{pr_th3}
By Proposition~\ref{prop4}, $R_n(\mathbb{P}) \le \sup_{x,y \in \mathbb{R}_{++}}\phi^{\mathbb{P}}(x,y)$.
Next, we prove the opposite inequality. Consider $\mathbb{P} \in \mathcal{C}_n$. We start with the case $n=2$ and  proceed similarly to Lu and Yu~\cite{LuYu2}. Denote the bivariate marginal distribution of $\mathbb{P}$ by $\mathbb{P}_2$. Consider a time matrix $T \in \mathbb{R}_{++}^{2\times 2}$ and denote $\frac{T_{11}}{T_{21}}$ by $x$,  $\frac{T_{12}}{T_{22}}$ by $y$. Construct the following matrix $T_0$:  \\
\begin{center}
$\begin{matrix}
 & \text{Task 1} & \text{Task 2}  \\
\text{Machine 1} \ &1& y \\
 \text{Machine 2} \ &\sfrac{1}{x} & 1
\end{matrix}$
\end{center}
The expected makespan of $\mathcal{A}^{\mathbb{P}_2}$ on this instance is $M(\mathbb{P}_2,T_0)$:
\begin{align*}
M(\mathbb{P}_2,T_0)   =  & \hspace{ 0.1 cm}  \mathbf{P}_{\mathbf{z}\sim \mathbb{P}_2}\big[z_1 > x, \ z_2 > y\big ](1+y)+\mathbf{P}_{\mathbf{z}\sim \mathbb{P}_2}\big [z_1 > x, \ z_2 \le y\big ] \\
&+\mathbf{P}_{\mathbf{z}\sim \mathbb{P}_2}\big [z_1 \le x, \ z_2 > y\big ] \max\big \{\sfrac{1}{x}, \ y \big \} \\
&+\mathbf{P}_{\mathbf{z}\sim \mathbb{P}_2}\big [z_1 \le x, \ z_2 \le y\big ]\big(1+\sfrac{1}{x}\big)\\
=&\hspace{ 0.1 cm} \big [1-H(x,y)-\big (F(y)-H(x,y)\big )-\big(F(x)-H(x,y)\big ) \big ](1+y)\hspace{ 1 cm} \\
&+\big (F(y)-H(x,y)\big )+\big (F(x)-H(x,y)\big )\max\big \{\sfrac{1}{x},  y \big \}\\
&+H(x,y)\big(1+\sfrac{1}{x}\big ) \\
=&\ 1+y-F(x)\big (1+y-\max\big \{\sfrac{1}{x}, \ y \big \}\big )-yF(y)\\
&+H(x,y)\big (y-\max\big \{\sfrac{1}{x}, \ y \big \}+1+\sfrac{1}{x}\big )\\
=& \ 1+y-\min\big\{ 1, \ 1-\sfrac{1}{x}+y\big\}F(x)-yF(y)\\
& +\min\big\{1+\sfrac{1}{x},  1+y\big\}H(x,y)\\
=& \ \phi^{\mathbb{P}}(x,y).
\end{align*}
Denote the minimum makespan on $T_0$ by $M^*$. By construction $M^* \le 1$, hence
\begin{align*}
&R_2(\mathbb{P}_2) \ge  R_2(\mathbb{P}_2,T_0)=\frac{M(\mathbb{P}_2,T_0)}{M^*}\ge  M(\mathbb{P}_2,T_0)=\phi^{\mathbb{P}}(x, y).
\end{align*}
This holds for all $x, y \in \mathbb{R}_{++}$, therefore
\begin{align*}
&R_2(\mathbb{P}_2) \ge  \sup_{x,y \in \mathbb{R}_{++}}\phi^{\mathbb{P}}(x,y).
\end{align*}
Now, fix $n>2$. Choose a small $\epsilon>0$ and consider the following time matrix $T_\epsilon$:
\begin{center}
$\begin{matrix}
 & \text{Task 1} & \text{Task 2} & \text{Task 3} &\dots& \text{Task n}\\
\text{Machine 1} \ &1& y & \epsilon &\dots& \epsilon \\
 \text{Machine 2} \ &\sfrac{1}{x} & 1 & \epsilon &\dots& \epsilon
\end{matrix}$
\end{center}
The expected makespan of $\mathcal{A}^{\mathbb{P}}$ on this instance, $M(\mathbb{P},T_\epsilon)$, satisfies
$$M(\mathbb{P}_2,T_0) \le M(\mathbb{P},T_\epsilon) \le M(\mathbb{P}_2,T_0)+(n-2)\epsilon,$$
and the optimal makespan, $M^*(T_\epsilon)$, satisfies
$$M^* \le M^*(T_\epsilon) \le M^*+(n-2)\epsilon.$$
Using the result for $n=2$,
\begin{align*}
&R_n(\mathbb{P})\ge\lim_{\epsilon \rightarrow 0}  R_n(\mathbb{P},T_\epsilon)=\lim_{\epsilon \rightarrow 0}  \frac{M(\mathbb{P},T_\epsilon)}{M^*(T_\epsilon)}=\frac{M(\mathbb{P}_2,T_0)}{M^*}\ge\phi^{\mathbb{P}}(x,y),
\end{align*}
which holds for all $x, y \in \mathbb{R}_{++}$.



\subsection{Possible critical points computation for the function \eqref{phi_f}} \label{ap1}

In this section we find possible critical points for the function
\begin{align*}
\phi^F(x,y)= y-\sfrac{1}{x} +\big(1+\sfrac{1}{x}-y\big)F(y)+\sfrac{1}{x}F(x),
\end{align*}
with $F$ defined in \eqref{f0} using \eqref{f_def} on $I_i \times I_j$ such that $i,j \in [2k], \ i+j \ge 2k+1$. By construction, in each interval $F$ is differentiable, and $\phi^F$ is differentiable on $I_i \times I_j$ (the expression simplifies for $i=1$ and $j=2k+1$). The first derivatives of the function $\phi^F(x,y)$ are:
\begin{align*}
\frac{\partial \phi^F(x,y)}{\partial x}= &\frac{1}{x^2} \bigg(1-F(x)-F(y)+\frac{\partial F(x)}{\partial x} x\bigg)\ \text{ and } \\
\frac{\partial \phi^F(x,y)}{\partial y}= & 1-F(y)+\big(1+\sfrac{1}{x}-y\big)\frac{\partial F(y)}{\partial y}.
\end{align*}
Next, we consider the three possible cases for $i,j$. For each of these cases we substitute $F(x), \ F(y)$ from \eqref{f0} and find the analytical solution to the system $\frac{\partial \phi^F(x,y)}{\partial x}=0, \ \frac{\partial \phi^F(x,y)}{\partial y}=0$. For this purpose we use Wolfram$|$Alpha~\cite{WA}. The obtained solution is denoted by $(x^*, \ y^*)$. \\

\noindent \textbf{\textit{Case 1.}} $k<i\le 2k, \ 1 \le j\le k$.\\

\noindent In this case $F(x)=1-{c^0_i}-{c^1_i}x, \ F(y)= {c^0_j}+\sfrac{{c^1_j}}{y}$ and $y \in (0,1)$. Hence
\begin{align*}
\frac{\partial \phi^F(x,y)}{\partial x}= &\frac{1}{x^2} \bigg({c^0_i}-{c^0_j}-\frac{{c^1_j}}{y}\bigg), \  \ \frac{\partial \phi^F(x,y)}{\partial y} \ge 0.
\end{align*}
The latter holds since $F(y) \le 1$, $F(y)$ is non-decreasing by construction, and $y \in (0,1)$. The sign of the derivative with respect to $x$ does not depend on $x$. The function $\phi^F(x,y)$ is non-decreasing in $y$ and is either non-increasing or non-decreasing in $x$. We not know this in advance, so we use the set $\big  \{(s_{i-1},s_j),(s_{i},s_j)\big  \} $ as possible critical points. \\

\noindent \textbf{\textit{Case 2.}} $k<i\le 2k, \ k<j\le 2k$.\\

\noindent In this case $F(x)=1-{c^0_i}-{c^1_i}x, \ F(y)= 1-{c^0_j}-{c^1_j}y$ and $y\ge 1$. Hence
\begin{align*}
\frac{\partial \phi^F(x,y)}{\partial x}= &\sfrac{1}{x^2} \big({c^0_i}-1+c^0_j+{c^1_j}y \big),\ \
\frac{\partial \phi^F(x,y)}{\partial y}= c^0_j-\big(1+\sfrac{1}{x}-2y \big){c^1_j}.
\end{align*}
As in \textbf{\emph{Case 1}}, the sign of the derivative with respect to $x$ does not depend on $x$, hence $\phi^F(x,y)$ is non-increasing or non-decreasing in $x$. We do not know this in advance, so we start with the set $\{s_{i-1},s_{i}\} \times \{s_{j-1},s_{j},y^*\}$ as possible critical points. If $i+j=2k+1$ (that is, the line $xy=1$ crosses the rectangle $I_i\times I_j$), we additionally consider the point $(\tfrac{1}{y^*}, y^*)$. We check all pairs for feasibility and exclude the infeasible ones. \\

\noindent \textbf{\textit{Case 3.}} $1 \le i\le k, \ k<j\le 2k$.\\

\noindent In this case $F(x)=c^0_i+\sfrac{{c^1_i}}{x}, \ F(y)=1-{c^0_j}-{c^1_j}y$, and
\begin{align*}
\frac{\partial \phi^F(x,y)}{\partial x}= &\frac{1}{x^2} \bigg(-c^0_i+{c^0_j}+{c^1_j}y - \frac{2c^1_i}{x} \bigg),\ \
\frac{\partial \phi^F(x,y)}{\partial y}= c^0_j-\big(1+\sfrac{1}{x}-2y\big)c^1_j.
\end{align*}
The sign of the derivatives is unknown, so we start with the set $\{s_{i-1},s_{i},x^*\} \times \{s_{j-1},s_{j},y^*\}$ as possible critical points. If $i+j=2k+1$ (that is, the line $xy=1$ crosses the rectangle $I_i\times I_j$), we additionally consider the set $\big \{(x^*,\tfrac{1}{x^*}),(\tfrac{1}{y^*}, y^*)\big \}$. We check all resulting pairs for feasibility and exclude the infeasible ones. \\

When $x \in I_1$ or $y \in I_{2k}$,  $\phi^F(x,y)$ simplifies by construction of \eqref{f0}. In our computations we analyze these situations separately.\\\\


\noindent \textbf{Acknowledgements}
We would like to thank the two anonymous reviewers for pointing out an inconsistency in the first proof of Corollary~\ref{cor:limit} and suggesting numerous improvements in the presentation of the results in this paper. This is a pre-print of an article published in Theory of Computing Systems. The final authenticated version is available online at: \href{https://doi.org/10.1007/s00224-019-09927-x}{https://doi.org/10.1007/s00224-019-09927-x}.

\bibliographystyle{spmpsci}
\bibliography{references_scheduling2} 



\end{document}